\documentclass[11pt,a4paper,reqno]{amsart}
\usepackage{amsmath,amsfonts,pictex,mathptmx,hyperref}

\newcommand{\re}{\mathbb{R}} \newcommand{\N}{\mathbb{N}}              
\newcommand{\zz}{\mathbb{Z}} \newcommand{\C}{\mathbb{C}}      

\newcommand{\Z}{{\zz}^n}                
\newcommand{\zn}{\Z}
\newcommand{\n}{{\N}_0}                 
\newcommand{\no}{{\N}_0}                 
\newcommand{\non}{\no ^n}

\newcommand{\rn}{{\re}^n}
\newcommand{\rl}{{\re}^{n-1}}
\newcommand{\Rn}{{\re}^n}
\newcommand{\Rl}{{\re}^{n-1}}

\newcommand{\cs}{{\mathcal S}}              
\newcommand{\cd}{{\mathcal D}}              

\newcommand{\bspq}{B^{s,a}_{p,q}\,}           
\newcommand{\fspq}{F^{s,a}_{p,q}\,}           
\newcommand{\bspqa}{B^{s,a}_{p,q}\,}           
\newcommand{\fspqa}{F^{s,a}_{p,q}\,}           
\newcommand{\aspq}{{A}^{s,a}_{p,q}\,}           

\newcommand{\np}{\frac{n}{p}}

\newcommand{\ep}{\frac{1}{p}}
\newcommand{\eq}{\frac{1}{q}}

\newcommand{\supp}{\operatorname{supp}}

\newcommand{\op}[1]{\operatorname{#1}}

\newcommand{\ffi}{\varphi}
\newcommand{\ffij}{\varphi _j}

\newcommand{\jsum}{\sum\limits _{j=0}^{\infty}}

\newcommand{\niusum}{\sum\limits _{\nu =0}^{\infty}}

\newcommand{\msum}{\sum\limits _{m\in \zn}}

\newcommand{\qnium}{Q^{a}_{\nu m}}

\newcommand{\imb}{\hookrightarrow}

\newcommand{\cf}{{\mathcal F}}
\newcommand{\cfi}{{\cf}^{-1}}

\newcommand{\be}{\begin{equation}}
\newcommand{\ee}{\end{equation}}
\newcommand{\beq}{\begin{eqnarray}}
\newcommand{\eeq}{\end{eqnarray}}
\newcommand{\beqq}{\begin{eqnarray*}}
\newcommand{\eeqq}{\end{eqnarray*}}

\newtheorem{satz}{Theorem}
\newtheorem{MainThm}{Main Theorem} 
\newtheorem{defi}{Definition}
\newtheorem{cor}{Corollary}
\newtheorem{lem}{Lemma}
\newtheorem{prop}{Proposition}
\newcounter{rem}[section] \renewcommand{\therem}{\thesection.\arabic{rem}}
\newenvironment{rem}{\par\medskip\noindent\refstepcounter{rem}%
{\bf Remark~{\therem}.}}{\par\medskip}

\catcode`\"=\active \let"\"

\renewcommand{\a}{|a|}
\renewcommand{\aa}{|a'|}
\newcommand{\ap}{\frac{|a|}p}
\newcommand{\anp}{\frac{a_n}p}

\newcommand{\g}{\gamma_{@!@!@!@!@!@!0}}

\renewcommand{\g}{\gamma_{\hspace{-0.1pt}0} }

\title[Traces\,---\,a complete treatment of the borderline
cases]{Traces of anisotropic Besov--Lizorkin--Triebel spaces\\
\,---\,a complete treatment of the borderline cases}
\author[W.~Farkas]{Walter Farkas$^*$}
\address{\small
 Inst.\ theor.\ Informatics and Mathematics,
 Univ.\ Fed.\ Armed Forces Munich,
 Werner-Heisenberg-Weg 39,
 D-85577 Neubiberg, Germany}
\author[J.~Johnsen]{Jon Johnsen}
\address{\small Department of Mathematical Sciences,
 Aalborg University,
 Fredrik~Bajers Vej 7E,
 DK--9220 Aalborg O, Denmark }
\author[W.~Sickel]{Winfried Sickel}
\address{\small Mathematics Department,
 Friedrich-Schiller-University Jena,
 Ernst-Abbe-Platz 1--4,
 D--07743 Jena, Germany} 

\thanks{$^*$Partly supported by the {\em Graduiertenkolleg\/} 
``Analytic and Stochastic Structures and Systems'' at FSU Jena, partly by
the DFG-project Ja~522/7-1.\\[4\jot]
{\tt Appeared in Mathematica Bohemica, vol.\ {\bf 125} (2000), no.~1, 1--37.}%
}


\begin{document} 

\begin{abstract} Including the previously untreated borderline cases, the
trace spaces (in the distributional sense) of the Besov--Lizorkin--Triebel
spaces are determined for the 
anisotropic (or quasi-homogeneous) version of these classes. The ranges of
the trace are  in all cases shown to be approximation spaces, and these
are shown to be different from the usual spaces precisely in the previously
untreated cases.
To analyse the new spaces, we
carry over some real interpolation results as well as the refined Sobolev
embeddings of J.~Franke and B.~Jawerth to the anisotropic scales.
\end{abstract}

\maketitle

\section{Introduction}
\enlargethispage{2\baselineskip}\thispagestyle{empty}

In this paper we present a complete solution of the trace problem
for the anisotropic (or rather quasi-homogeneous) Besov and Lizorkin--Triebel
spaces, denoted by $\bspq$ and $\fspq$, respectively. The definitions are
recalled in Appendix~\ref{spaces-app} below.

Here the trace is the operator $\g$,
\begin{equation}\label{1-1}
  f(x_1,x_2, \ldots , x_{n-1},x_n) 
  \mathrel{\rule[-.06pt]{.5pt}{5pt}%
                               \hspace{-0.7pt}\mathalpha{\xrightarrow{\mbox{$\ \g\:\;$}}}}
  f(x_1,x_2, \ldots ,x_{n-1},0) ,
\end{equation}
which restricts functions on $\rn$ to the hyperplane $\Gamma = \{ x_n =0\}$ in
$\rn$ ($n \ge 2$)\,---\,in general this is defined in the obvious way on the subspace
$C(\re,\cd'(\rl))$ of $\cd'(\rn)$. For the $\bspq(\rn)$ and $\fspq(\rn)$
under consideration,
this coincides with the extension by continuity from the Schwartz space
$\cs(\rn)$, except for $p=\infty$ and $q=\infty$ in which cases $\cs$ is not dense;
however, because of embeddings, the latter exception is only felt for
$F^{a_n,a}_{1,\infty}$, where $s=a_n$ is the lowest possible value, $a_n$
being the modulus of anisotropy associated to $x_n$. (For simplicity, $\Rn$ is often suppressed in
$\bspq$ and $\fspq$ etc.\ when confusion is unlikely to result.)

The trace problem consists in finding spaces $X$ and $Y$, as subspaces of
$\cd'(\rn)$ and $\cd'(\rl)$, respectively, such that $\g$
yields a continuous, linear \emph{surjection}
\begin{equation}
  \g\colon X\to Y.
\end{equation}
In this paper we determine those $\bspq$ and $\fspq$ which allow such a
$Y$ to be found, and we moreover determine the optimal $Y$ for these
choices of $X$. (The existence of a right inverse of $\g$ is also discussed.) 

One effect of allowing \emph{anisotropic} spaces is that $\g$ is
studied on larger domains. However, the main motivation for the anisotropic spaces is
that they are indispensable for the fine theory of parabolic boundary
problems. For example it is well known that in a treatment of
$\partial_t-\Delta$ it is necessary 
to use $\bspq$ and $\fspq$ for $a=(1,\dots,1,2)$ 
(gathering the moduli of anisotropy to form a vector in the $(x,t)$-space) and 
that $F^{s,a}_{p,2}$ (locally) equals the intersection of
$L_p(\re,H^s_p(\rl))$ and $H^{s/2}_p(\re,L_p(\rl))$. 
We refer the reader to works of G.~Grubb \cite{G1,G2} for a recent
treatment, based on L.~Boutet de~Monvel's pseudo-differential calculus, of
parabolic initial boundary problems in anisotropic 
Besov and Bessel potential spaces.

If desired, the reader may 
specialise to the isotropic case, which is given by $a=(1,\dots,1)$.

\bigskip

In the following review of results, comparison with other
works is often postponed (for simplicity's sake) to
later sections.

\begin{figure}[htbp]
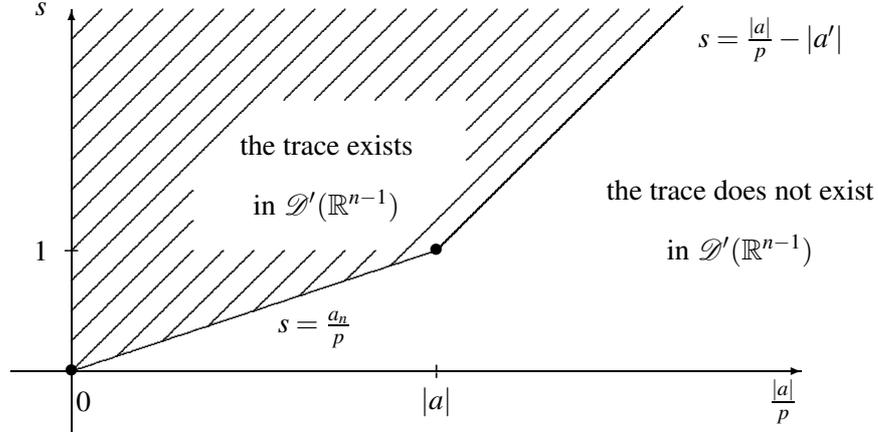
 
\hfill
\beginpicture

\setcoordinatesystem units <0.8cm,0.8cm>
\unitlength0.8cm

\put {\vector (0,1) {7} } [Bl] at 1 0
\put {\vector (1,0) {13} } [Bl] at 0 1

\plot 1 1 7 3 /
\plot 7 3  11 7 /
\plot 7 0.9 7 1.1 /
\plot 0.9 3 1.1 3 /

\put {$s$} at 0.5 7

\put {$\frac{\a}{p}$} at 12.7 0.5 
\put {${\a}$} at 7 0.5 
\put {${1}$} at 0.5 3 
\put {${0}$} at 1.2 0.5 
\put {$s= \frac{\a}{p} - \aa$} at 12.5 6.5 
\put {$s= \frac{a_n}{p} $} at 5 1.7 
\put {$\bullet$} at 1 1
\put {$\bullet$} at 7 3 
\put {the trace exists} at 5.2 4.75 
\put {in $\cd'(\re^{n-1})$} at 5.2 3.75 
\put {the trace does not exist} at 12 4
\put {in $\cd'(\re^{n-1})$} at 12 3



\plot 1 2 3 4 /
\plot 1 3 5 7 /
\plot 1 4 4 7 /
\plot 1 5 3 7 /
\plot 1 6 2 7 /
\plot 1 2.5 5.5 7 /
\plot 1 3.5 4.5 7 /
\plot 1 4.5 3.5 7 /
\plot 1 5.5 2.5 7 /
\plot 1 6.5 1.5 7 /
\plot 1 1.5 3 3.5 /
\plot 1 1 3 3 /
\plot 4.5 5.5 6 7 /
\plot 5 5.5 6.5 7 /
\plot 5.5 5.5 7 7 /
\plot 6 5.5 7.5 7 /
\plot 6.5 5.5 8 7 /
\plot 7 5.5 8.5 7 /
\plot 7.5 5.5 9 7 /
\plot 7.5 5 9.5 7 /
\plot 7.5 4.5 10 7 /

\plot 1.75 1.25 3.5 3 / 
\plot 2.5 1.5 4 3 /
\plot 3.25 1.75 4.5 3 /
\plot 4 2 5 3 /
\plot 4.75 2.25 5.5 3 /
\plot 5.5 2.5 6 3 /
\plot 6.25 2.75 10.5 7 /

\thicklines
\plot 7 3 11 7 /
\plot 7.05 3.05 11.05 7.05 /

\endpicture
\hfill
\caption{The region specified by formula \eqref{1-2}}
  \label{param-fig}
\end{figure}

For the continuity of $\g$ from $\bspq$ or $\fspq$ to
$\cd'(\rl)$ two different conditions are necessary for $p\ge1$
and $p<1$. Introducing $\a:=a_1+\dots+a_n$ and similarly 
$\aa=a_1+\dots+a_{n-1}$ for the spaces over $\rl$ (so that in the
isotropic case $\a=n$ and $\aa=n-1$), these may be expressed in the following
way by letting $t_+:=\max(0,t)$:
\begin{equation}   \label{1-2}
s\ge \max\bigl(\anp,\ap-\aa\bigr)= \anp +\aa\bigl(\ep-1\bigr)_+;
\end{equation}
in other words the correction $\aa(\ep-1)$ appears for $0<p<1$.
In the case of equality it is moreover necessary that $q\le1$ in the $B$-case
and $p\le1$ in the $F$-case; cf.~Figure~\ref{param-fig}.

This has earlier been known to specialists (we carry over explicit
isotropic counter-examples of the second author \cite[Rem.~2.9]{Jo2}). The
case $s=\anp$ for $p\ge1$ was investigated by 
V.~I.~Burenkov and M.~L.~Gol'dman \cite{BG}, however only for $q=1$;
in the isotropic case $0<q\le1$ and, for $p<1$, the borderline case
$s=\np-(n-1)$ was treated by J.~Johnsen~\cite{Jo3}, but 
surjectivity  was left open for $p<1$. The present article may therefore be
seen as a continuation of \cite{BG} and \cite{Jo3}. Emphasis will be on the
borderline cases mentioned, since it is known (and
comparatively easy) that $\g$ is a bounded, right invertible surjection
\begin{equation}
  \g\colon \bspq(\Rn) \to B^{s-\anp,a'}_{p,q}(\Rl),\qquad
  \g\colon \fspq(\Rn) \to B^{s-\anp,a'}_{p,p}(\Rl)
\end{equation}
in any of the \emph{generic} cases (i.e.\
those with strict inequality in \eqref{1-2}). This is, of course, the
well-known loss of $\ep$ in the isotropic case.

For $s=\ap-\aa$ with $p<1$ the optimal $Y$'s are determined below, and it
turns out that whenever $p$, $q\le1$ they are neither Besov nor Lizorkin--Triebel
spaces. So  in order to complete the range characterisation, we introduce
another scale $\aspq$, which previously has been 
investigated mainly by Russian specialists in function spaces. 
 
In fact, $\g$ still lowers $s$ by $\anp$ and it is a bounded surjection
\begin{align}
  \g\colon B^{\ap-\aa,a}_{p,q}(\Rn)\to & A^{\aa(\ep-1),a'}_{p,q}(\Rl) 
  \quad\text{for}\quad p,q\in\,]0,1], 
  \label{Blim-eq} \\
  \g\colon F^{\ap-\aa,a}_{p,q}(\Rn)\to & A^{\aa(\ep-1),a'}_{p,p}(\Rl) 
  \quad\text{for}\quad 0<p\le1,\quad 0<q\le\infty. 
  \label{Flim-eq}
\end{align}
To have instead e.g.\ a Besov space as the codomain, one can use an embedding of
$A^{\aa(\ep-1),a'}_{p,q}(\Rl)$ into 
$B^{\aa(\frac1r-r),a'}_{r,\infty}(\Rl)$ for $r=\max(p,q)$ (which is optimal).
An investigation of the borderline cases is given in
Section~\ref{collect-ssect} below.

However, because of the various identifications between the $\aspq$ and
the Lebesgue, Besov and Lizorkin--Triebel spaces etc., it is possible to
formulate \emph{all} trace results in a concise way in terms of the $\aspq$:

\begin{MainThm}   \label{main-thm}
For a given anisotropy $a=(a_1,\dots,a_n)$, the following assertions are equivalent: 
\begin{itemize}
  \item[\textup{(a)}]
    the operator $\g$ is a continuous mapping from $\bspq(\Rn)$ into 
    $\cd'(\rl)$;
  \item[\textup{(b)}]
    the operator $\g$  maps  $\bspq(\Rn)$ continuously {\bf onto\/}
    $A^{s-\frac{a_n}{p},a'}_{p,q}(\Rl)$;
  \item[\textup{(c)}] the triple $(s,p,q)$ satisfies $s\ge
    \anp+\aa(\ep-1)_+$ and, in case of equality, also $0<q\le 1$.
\end{itemize}
For the Lizorkin--Triebel spaces $\fspq$ the analogous result holds if
one replaces $A^{s-\frac{a_n}{p},a'}_{p,q}$
by $A^{s-\anp,a'}_{p,p}$ in \textup{(b)} and replaces `$q\le1$'
by `$p\le1$' in \textup{(c)}. 
\end{MainThm}

It should be noted that the formal introduction of $\aspq$ in
Definition~\ref{aspq-defn} below allows us to give a short, self-contained
proof of the main theorem in the $\bspq$-case (thus a unified proof of all
underlying borderline cases); the Lizorkin--Triebel case is then deduced
from the Besov case by establishing a certain $q$-independence. This partly follows
the work of M.~Frazier and B.~Jawerth \cite{FJ2}, but we point out
and correct a flaw in the proof of \cite[Th.~11.1]{FJ2}, see
Remark~\ref{FJ-rem} below.

In order to deduce the relations between $A^{\aa(\ep-1)}_{p,q}(\Rl)$,
which enter in \eqref{Blim-eq}--\eqref{Flim-eq} above, and the usual spaces,
we need anisotropic versions of the optimal mixed Sobolev embeddings between
the two scales $\bspq$ and $\fspq$.

For this purpose we carry over these embeddings (due to B.~Jawerth and
J.~Franke), hence also some necessary real interpolation results for
$\fspq$, to the anisotropic setting. See Appendix~\ref{atom-app} below for these results.

\bigskip

In Section~\ref{gen-sect} below we introduce a working definition (based on
a limit) of $\g$ and then present results for the generic cases. The
borderline cases are treated in Section~\ref{bdr-sect}, in particular the range
spaces are presented for the cases with $s=\ap-\aa$ in terms of the
approximation spaces $\aspq$, which are formally introduced there for this
purpose. The relations between $\aspq$ and the Besov and Lizorkin--Triebel
spaces are elucidated in Section~\ref{apr-sect}, and the proofs of the assertions
are to be found in Section~\ref{proof-sect}. The appendices collect the notation and the
necessary facts about Besov and Lizorkin--Triebel spaces, in particular the
extension to the anisotropic case of some well-known facts.

\begin{rem}
We should emphasise that, for $p<1$, also a different
operator has been studied under the label `trace'.
Indeed, when $\g$ is restricted to the Schwartz space
$\cs(\rn)$, there are extensions by continuity $T\colon \bspq\to Y$ at least
if $s>\ep$; while this is effectively weaker than \eqref{1-2}, it may only
be obtained by taking $Y$ outside of $\cd'$, and it was shown in \cite{Jo3}
that $T$ is different from $\g$ since
\begin{equation}
  \g(\varphi(x')\otimes\delta_0(x_n))=\varphi,\qquad
         T(\varphi(x')\otimes\delta_0(x_n))=0     
\end{equation}
for all $\varphi\in \cs(\rl)$.

Moreover, using the sharp result for elliptic boundary operators in Besov
and Lizorkin--Triebel spaces obtained in \cite{Jo2}, it was also
proved in \cite{Jo3} that $T$ is unsuitable for the study 
of elliptic boundary problems. We shall therefore not consider this
other possibility here; it was discussed at length by M.~Frazier and
B.~Jawerth~\cite{FJ1, FJ2} and H.~Triebel~\cite{Tr6, Tr9}.
\end{rem}
\begin{rem}   \label{param-rem}
It is noted once and for all that we consider arbitrary $s\in\re$ and
$p$, $q\in\,]0,\infty]$, although $p<\infty$ is to be understood throughout
for the $\fspq$ spaces; and all such admissible parameters are 
considered unless further restrictions are stated.
\end{rem}

\section{Generic properties of the trace}
  \label{gen-sect}
To set the scene properly, we introduce the trace in a formal way.
The reader should consult Appendix~\ref{notation-app}--\ref{spaces-app}
first for the anisotropic spaces and for the corresponding anisotropic
distance $|x|_a$ and dilation $t^ax=(t^{a_1}x_1,\dots,t^{a_n}x_n)$, both
defined on $\rn$; here and throughout $a=(a_1,...,a_{n-1},a_n)=(a',a_n)$ will be a
given anisotropy. 

In addition, $\cf$ denotes the Fourier transform and $\cfi$ the inverse,
extended from the Schwartz space $\cs(\rn)$ to its dual $\cs'(\rn)$. 
In different dimensions, say  $\rl$ and $\re$, Fourier transformation will be
indicated by $\cf_{n-1}$ and $\cf_1$, respectively. 
Let $\psi \in C_0^\infty (\rn )$ be a function such that
\begin{equation}
  \label{2-0}
  \psi (x) =1 \qquad \mbox{if} \quad |x|_a\le 1 \qquad \mbox{and} \qquad
  \psi (x)=0 \quad \mbox{if} \quad |x|_a \ge 2.
\end{equation}
For such $\psi$, we may define a smooth, anisotropic
dyadic partition of unity $(\ffij )_{j\in \no}$ by letting $\ffi_0(x)=\psi(x)$ and
\begin{gather}   
  \ffij (x) = \ffi _{0}(2^{-ja}x)-
  \ffi _{0}(2^{(-j+1)a}x)\quad\mbox{if}\quad j\in \N .
  \label{22-2}
\end{gather}
Indeed, since $\psi(2^{Na}\xi)=\sum_{j=0}^N \varphi_{j}(\xi)$, it is clear that
\begin{equation}
 \label{22-3}
   \jsum \ffij (x)=1 \quad \mbox{for}\quad x\in \rn.
\end{equation}

\subsection{The working definition}

Using a fixed $\psi$ of the above type, we have for all $f \in \cs' (\rn)$ 
\begin{equation}
  \label{2--1}
  f = \sum_{j=0}^{\infty} \cfi [\ffij\cf f ] 
  \qquad (\text{convergence in $\cs'$}) .
\end{equation}
Since,  by the Paley--Wiener--Schwartz Theorem, $\cfi [\ffij\cf f ]$ is
continuous,  the trace has an immediate
meaning for this function. As a temporary working definition we
therefore let
\begin{equation}\label{2-01}
  \g f = \lim_{N\to \infty} \sum_{j=0}^N \g (\cfi [\ffij\cf f] ) 
\end{equation}
whenever this limit exists in $\cs'(\rl)$.  However, we should make the
following remarks to this definition.

On the one hand, it is possible to show the next result (which summarises the
most well-known facts on $\g$) by relatively simple arguments:

\begin{satz}
  \label{gen-thm}
The operator $\g$ maps $\bspq (\Rn)$ (or $\fspq$) continuously into
$\cs'(\rl)$ only if $(s,p,q)$ satisfies either
\begin{equation}\label{2-1}
  s > \tfrac{a_n}{p} + |a'| (\tfrac{1}{p} - 1 )_+
\end{equation}
or, alternatively,
\begin{equation}
  \label{2-2}
  s = \tfrac{a_n}{p} + |a'| (\frac{1}{p} - 1 )_+ 
  \qquad \text{and} \qquad 0< q \le 1 
  \quad\text{(respectively $p\le 1$)} . 
\end{equation}
When \eqref{2-1} holds, then $\g$ is actually a continuous map
$\bspq(\Rn)\to B^{s-\anp,a'}_{p,q}(\Rl)$ and $\fspq(\Rn)\to B^{s-\anp,a'}_{p,p}(\Rl)$.
\end{satz}

On the other hand, by the same line of thought as in \cite{Jo3}, one may, as
we show in this paper, deduce from the proof of Theorem~\ref{gen-thm} that
the just defined operator $\g$ coincides with (a restriction of) the map
\begin{equation}
  \label{r0-eq}
  r_0f:= f(0)\quad\text{for $f(t)$ in $C(\re,\cd'(\rl))$}.
\end{equation}
To conclude this relation between $\g$ and $r_0$, we apply the next
result where $C_{\op{b}}(\re, X)$ stands for the (supremum normed) space of uniformly
continuous, bounded functions valued in the Banach space $X$:

\begin{prop}
  \label{Cembd-prop}
When $\bspq(\Rn)$ and $\fspq(\Rn)$ have parameters $s$, $p$ and $q$
satisfying \eqref{2-1}, or the pertinent version of \eqref{2-2}, then
\begin{equation}
  \bspq(\Rn),\ \fspq(\Rn)\hookrightarrow C_{\op{b}}(\re,L_{p_1}(\rl))
  \label{Cb-incl}
\end{equation}
hold with $p_1=\max(p,1)$.
\end{prop}
Given this result (see Section~\ref{Cembd-pf} for the proof), it follows
from the fact that uniform convergence 
implies pointwise convergence (hence that $r_0$ is continuous from
$C_{\op{b}}(\re,L_{p_1}(\rl))$) that, for any $f$ in one of the spaces on
the left hand side of \eqref{Cb-incl},
\begin{equation}
  \label{Cb-eq}
  r_0f=\lim_{N\to\infty} r_0\sum_{j=0}^N \cfi(\ffij\cf f)=
  \sum_{j=0}^{\infty} \cfi(\ffij\cf f)(\cdot,0)=\g f.
\end{equation}
Indeed, for $q<\infty$ the decomposition in \eqref{2--1} converges in the
topology of the $\bspq$ or $\fspq$ space that contains $f$, hence there is
also convergence in $C_{\op{b}}(\re,L_{p_1}(\rl))$; with one exception it is
always possible to reduce to the case with $q<\infty$ by means of
embeddings, e.g.\ a space with $p\le 1$ is a subspace of $F^{a_n,a}_{1,q}$
for some $q<\infty$ according to \eqref{2-2}.

The just mentioned exception is the space $F^{a_n,a}_{1,\infty}(\Rn)$ for
which \eqref{Cb-eq} requires a sharper argument: for arbitrary $f\in
C_{\op{b}}(\re,L_{1}(\rl))$ one can take any $\eta\in\cs(\rn)$ with
$\int_{\rn}\eta\,dx=1$, let
$\eta_k(x):=k^{\a}\eta(k^ax)$ and then show that
\begin{equation}
  \label{eta-eq}
  \eta_k * f (x',0)\to f(x',0)\quad\text{in $L_1(\rl)$ for $k\to\infty$}.
\end{equation}
This is sufficient because it applies to any $f\in
F^{a_n,a}_{1,\infty}$ by Proposition~\ref{Cembd-prop}, and while
$r_0f=f(\cdot,0)$ by definition, $\sum_{j=0}^N \cfi(\ffij\cf f)(x)= \eta_k*
f(x)$ holds for $\eta=\cfi\psi$ and $k=2^N$, so that altogether the first
equality sign of \eqref{Cb-eq} is justified. 

However, \eqref{eta-eq} may be
verified by the usual convolution techniques, for if translation by $y'$ in
$\rl$ is denoted by $\tau_{y'}$ and $\|\cdot| L_1(\rl)\|$ is replaced by
$\|\cdot\|_1$, the translation invariance gives that 
\begin{equation}
  \begin{split}
  \|\,\eta_k * f (\cdot,x_n)-f(\cdot,x_n)\,\|_1 \le &
  \int_{\rn} |\eta(z)|\cdot
  \|\,f(\cdot,x_n-k^{-a_n}z_n)-f(\cdot,x_n)\,\|_1 \,dz
\\
  & +
  \int_{\rn} |\eta(z)| \cdot
  \|\,(\tau_{k^{-a'}z'}-I)f(\cdot,x_n)\,\|_1   \,dz.
  \end{split}
\end{equation}
Setting $x_n=0$, one may for any $\varepsilon>0$ take $c>0$ such that
$(\tau_{y'}-I)f(\cdot,0)$ has $L_1$-norm less than $\varepsilon$ when
$|y'|_{a'}<c$ (since $|\cdot|$ and $|\cdot|_{a'}$ define the same
topology in $\rl$). In the second integral it thus remains to control the
region where $|z'|_{a'}\ge ck$, but a majorisation by $2\|\,f(\cdot,0)\,|
L_1\|$ shows that this contribution is less than $\varepsilon$ for all
$k$ eventually; by the continuity with respect to $x_n$, and a similar
splitting, the first integral is also $<\varepsilon$ eventually. 

Hence \eqref{eta-eq} holds (the real
achievement is the less elementary proof of Proposition~\ref{Cembd-prop}),
and thus \eqref{Cb-eq} is proved for all spaces treated in the present paper.
In particular, this means that $\g u$ is independent of the choice of $\psi$
(and of the anisotropy $a$).

Summing up the above discussion, we have proved

\begin{cor}
  When the operators $r_0$ and $\g$ are defined as in \eqref{r0-eq} and
  \eqref{2-01} above, then 
\begin{equation}
  r_0f=\g f
\end{equation}
holds for all functions $f$ in the spaces $\bspq$ and $\fspq$ fulfilling
\eqref{2-1} or \eqref{2-2}.
\end{cor}

\begin{rem}
Our working definition of $\g f$ has been used since the late 1970's; cf.\
\cite{Ja,B,Tr6}. Nevertheless, the consistency with the trace on $C(\re,\cd'(\rl))$ was,
to  our knowledge, first proved for the Besov spaces in \cite{Jo3}. By
Proposition~\ref{Cembd-prop} and \eqref{Cb-eq} above this consistency holds
for all the considered spaces; the consistency extends to the trace defined
on the entire Colombeau algebra ${\mathcal G}(\rn)$, which contains $\cd'(\rn)$, see
\cite[Prop.~11.1]{Ober}. 
\end{rem}

\subsection{Linear extension}
  \label{K-ssect}

Taking, as we may, $\eta _0$ and $\eta \in \cs (\re )$ such that
$\supp \eta _0\subset \,]-1,1[$ and $\supp \eta \subset \,]1,2[$ 
and normalised so that 
\begin{equation}
\label{22-4}
(\cf _1^{-1}\eta _0)(0)=(\cf _1^{-1}\eta )(0)=1,
\end{equation} 
we set
\begin{equation}
\label{22-5}
\eta _{j}(t)=\eta (2^{-ja_n}t) \quad \mbox{for any}\quad t\in \re\quad
\mbox{if}\quad j\ge 1 .
\end{equation}
For any $v\in \cs '(\rl )$ such that the following series converges,  one
can now define an extension to $\rn$ by means of the partition of unity
$(\ffij )_{j\in \no}$ in \eqref{22-2}:
\begin{equation}
\label{22-6}
K v(x',x_n)= \jsum 2^{-ja_n} \, (\cf _{1}^{-1}\eta _{j})(x_n)\,
\cf _{n-1}^{-1}[\ffij (\cdot,0)\, \cf _{n-1}v](x').
\end{equation}
Using a homogeneity argument, we have
$ 2^{-ja_n} \, (\cf _{1}^{-1}\eta _{j})(0)=1$ for any $j\in \no$,
so the termwise restriction to $x_n=0$ gives, with weak convergence in
$\cs'(\rl )$, 
\begin{equation}   \label{22-7}
  \sum_{j=0}^N 2^{-ja_n} \, (\cf _{1}^{-1}\eta _{j})(0)\,
  \cf _{n-1}^{-1}[\ffij (\cdot,0)\, \cf _{n-1}v]
  \underset{N\to\infty}{\longrightarrow} \jsum  \cf _{n-1}^{-1}[\ffij (\cdot,0)\, \cf _{n-1}v]
  =v.
\end{equation}
By the working definition of $\g$, this means that $\g Kv$ is defined for
such $v$, hence $\g K=I$; i.e.\ $K$ is a linear extension.
 
When it is understood that the convergence of the series \eqref{22-6} is part of the
assertion, one has
\begin{satz}  \label{gen1-thm}
\begin{itemize}
  \item[\rm{(i)}] The operator $K$ maps $B^{s-\frac{a_n}{p},a'}_{p,q}(\Rl )$
  continuously into $\bspq (\Rn )$.
  \item[\rm{(ii)}]  For $0<p<\infty$ the operator $K$ maps
  $B^{s-\frac{a_n}{p},a'}_{p,p}(\Rl )$ continuously into $\fspq (\Rn )$.
\end{itemize}
\end{satz}

Here there are no restrictions in $s$, that is, the assertions in Theorem 2 hold for
all $s \in \re$ (which is to be expected since $K$ is a Poisson operator).

Since the relation $\g K=I$ was found above, one has as a consequence the
next result.
\begin{satz}
  \label{gen2-thm}
Let $s > \frac{a_n}{p} + |a'| (\frac{1}{p} - 1 )_+ $.
Then the operator $\g$ maps $\bspq (\Rn )$ continuously
{\bf onto\/} $B^{s-\frac{a_n}{p},a'}_{p,q}(\Rl )$, and for $0<p<\infty$ it
maps $\fspq(\Rn)$ onto $B^{s-\frac{a_n}{p},a'}_{p,p}(\Rl )$; $K$ is in both cases a
linear right inverse of $\g$. 
\end{satz}

Although Theorems~\ref{gen1-thm}--\ref{gen2-thm} above are unsurprising
(indeed, well-known in the isotropic case), they deserve to be compared with
the borderline results in the next section.

\begin{rem}
The contents of the above theorems are known to a wide extent for the
classical parameters $1\le p,q\le\infty$; cf.\ the works of O.~V.~Besov,
V.~P.~Ilyin and S.~M.~Nikol'skij \cite{BIN}, S.~M.~Nikol'skij \cite{Nik},
V.~I.~Burenkov and M.~L.~Gol'dman \cite{BG} and G.~A.~Kalyabin \cite{Ka1}.
For the isotropic case we also refer to the works of J.~Bergh and
J.~L{\"o}fstr{\"o}m \cite{BL}, M.~Frazier and B.~Jawerth \cite{FJ1,FJ2} and
H.~Triebel \cite{B,Tr6,Tr9} as well as to the remarks in \cite{Jo3} and in
the present paper. The study of the trace problem for $0<p<1$ was
initiated by B.~Jawerth \cite{Ja,Ja2}, but the first to find the
borderline $s=\frac np-(n-1)$ for $0<p<1$ seems to be either B.~Jawerth or
J.~Peetre \cite[Rem.~2.3]{Pe1}. (Peetre \cite[Note~1]{Pe1} actually gives
credit to \cite{Ja2} for this, and vice versa in \cite[Rem.~2.2]{Ja2}.) 
\end{rem}

\begin{rem}
The borderline $s=\frac{1}p$ itself was found by S.~M.~Nikol'skij in 1951
when he proved the continuity of $\g\colon B^s_{p,\infty}(\rn)\to
B^{s-\ep+(n-1)(\frac1r-\ep)}_{r,\infty}(\rl)$ for $s>\ep$ and any $r>p$ (the
result was actually anisotropic and valid for the restriction to linear
submanifolds of codimension $m\ge1$); cf.~\cite{Nik51} and also
\cite[6.5]{Nik}. Traces of Sobolev spaces $W^1_p$ were studied first by
N.~Aronszajn~\cite{Ar} around 1954.
Later, around 1957, E.~Gagliardo \cite{Ga} considered the trace of
$W^1_1$, which constitutes an `extremal' case; cf.\ the vertex in Figure~1. 
However, the first explicit counterexamples for the borderline
$s=1/p$ seem to be put forward by G.~Grubb \cite{GG} (who stated they were due to
L.~H{\"o}rmander); these necessary conditions were then expanded in
\cite{Jo2}, and in Lemma~\ref{711-l4} ff.\ below these are 
supplemented to a set of necessary conditions for the anisotropic Besov and
Lizorkin--Triebel spaces; in view of this paper the conditions are also
\emph{sufficient} for a solution of the distributional trace problem for the
spaces considered.
\end{rem}

\section{The borderline cases}   \label{bdr-sect}
In all remaining cases where $\g$ is a
continuous operator into $\cs'(\rl)$ it turns out that $\g$ has properties
different from the generic ones in Theorem~\ref{gen-thm}.
Recall that it remains to investigate
\begin{equation}
s = \left\{\begin{array}{lll} 
\frac{a_n}{p} & \qquad & \mbox{if} \quad 1 \le p \le \infty, \\
\frac{|a|}{p} - |a'| & & \mbox{if} \quad 0 < p < 1;
\end{array}\right.
\end{equation}
this amounts to the following five cases, see Figure~1, of which only the subcase $q=1$ of
the first two has been completely covered in the literature hitherto (whilst
only the first case and the subcase $q<\infty$ of the fourth have been settled isotropically):
\begin{itemize}
  \item $B^{\frac{a_n}{p},a}_{p,q} (\Rn)$ with $1 \le p < \infty$ and 
    $0 < q \le 1$;
  \item $B^{0,a}_{\infty, q}(\Rn)$ for  $0 < q \le 1$;
  \item $B^{\frac{|a|}{p}-|a'|,a}_{p,q}(\Rn)$ for $ 0 < p < 1$ and $0 < q \le 1$;
  \item $F^{a_n, a}_{1, q}(\Rn)$ with  $ 0 < q \le \infty$;
  \item $F^{\frac{|a|}{p} - |a'|, a}_{p,q}(\Rn)$ for
    $0 <  p < 1$ and $0 < q \le \infty$.
\end{itemize}
However, as a preparation for these cases, some
preliminaries are dealt with in the next subsection.

\subsection{Approximation spaces and nonlinear extension}

To describe the trace classes  in the limit situations we introduce a
another class of spaces, actually a half-scale, related to the approximation
by entire analytic functions of exponential type.

\begin{defi}   \label{aspq-defn}
Let $p, q\in\,]0,\infty]$ and let $(s,p,q)$ fulfil one of the following two conditions:
\begin{align}
  s &> \a \bigl( \ep - 1\bigr)_+ ,
  \label{spq-ineq}  \\
  s &= \a \bigl( \ep - 1\bigr)_+ \qquad \text{and}\qquad 0 < q \le 1.
  \label{spq-eq}
\end{align}
Then we define the anisotropic \emph{approximation} space $\aspq$ to be the set
\begin{multline}
  \aspq (\Rn)  =   \bigg\{ f \in \cs' (\rn) \bigm|
  \forall j\in\n,\exists h_j \in L_p (\rn) \cap \cs' (\rn), 
\\
  \qquad\supp \cf h_j \subset \{\xi: \: |\xi|_a \le 2^j\},
  \qquad f \underset{\cs'}{=} \sum_{j=0}^\infty \, h_j ,
\\
  \big( \sum_{j=0}^\infty \, 2^{jsq}\, \|\, h_j\, |L_p (\rn)
  \|^q\big)^{1/q} < \infty 
  \bigg\} ;
\end{multline}
this is equipped with the quasi-norm
\begin{equation}
  \| \, f \, | \aspq (\Rn)\| = \inf \, \big( \sum_{j=0}^\infty 
  2^{jsq}\, \|\, h_j\, |L_p(\rn)\|^q\big)^{1/q},  
\end{equation}
where the infimum is taken over all admissible representations $f=\sum h_j$.
(If $q=\infty$, the $\ell_q$-norm should be replaced by the supremum
over $j$ in both instances above.)
\end{defi}

S.M.~Nikol'skij and O.V.~Besov have (in connection with Besov
spaces for $1\le p\le\infty$) defined such spaces 
$\aspq$. For a comprehensive treatment and additional references, see the
monograph of S.M.~Nikol'skij \cite[3.3 and 5.6]{Nik}; cf.\
also H.~Triebel \cite[2.5.3]{Tr6}.  

The restrictions on $s$ make the $\aspq (\rn)$ continuously embedded into
$\cs'(\rn)$, cf. the following proposition and its corollary on
the identifications between the $\aspq$ and other classes of functions.
Here and throughout $C_{\op{b}}(\rn)$ denotes the set of all uniformly
continuous, bounded, complex-valued functions on $\rn$ equipped with the
supremum norm. 

\begin{prop}   \label{aspq-prop}
\begin{itemize}
  \item[{\rm (i)}] 
    $\aspq (\Rn) = \bspq (\Rn)$ if, and
    only if,  $s >   |a| \, (\ep -1)_+$.
  \item[{\rm (ii)}] Let $1 \le p < \infty$ and $0 < q \le 1$.
    Then $A^{0,a}_{p,q} (\Rn)= L_p (\rn) $.
  \item[{\rm (iii)}] Let $0 < q \le 1$.
    Then $A^{0,a}_{\infty, q} (\Rn)= C_{\op{b}}(\rn) $.
\end{itemize}
\par\noindent
In the affirmative cases the quasi-norms are equivalent.
\end{prop}
\begin{rem}
Parts (ii) and (iii) are contained in Burenkov's and Gol'dman's work
\cite{BG}, at least implicitly. An isotropic counterpart is stated in Oswald
\cite{Os}, albeit with the spaces based on approximation by splines. 
Furthermore, the ``if''-part in (i) 
has been known before, cf. Nikol'skij \cite[3.3 and 5.6]{Nik},
Triebel~\cite[2.5]{Tr6}; recently also Netrusov~\cite{Ne} considered this
issue. 
\end{rem}

In virtue of Proposition~\ref{aspq-prop} the continuity of
$\aspq\hookrightarrow\cs'$ is clear for $p\ge1$; when $0<p<1$ the
anisotropic Nikol'skij--Plancherel--Polya inequality (cf.~\cite[2.13]{Ya})
and the restrictions \eqref{spq-ineq}--\eqref{spq-eq} immediately give an
embedding into $L_1$, so one has

\begin{cor}   \label{aspq-cor}
The classes $\aspq (\Rn)$ are continuously embedded into
$L_{\max(1,p)}(\rn)$ (which for any $p\in\,]0,\infty]$ is a subspace of
$\cs'$).
\end{cor}

For a more detailed comparison of the classes $\aspq (\Rn)$ in the
remaining limit situations with spaces of
Besov--Lizorkin--Triebel type we refer to Section 4.

\bigskip

To make extensions to $\rn$ of suitable $f$ in $\cs'(\rl)$, we consider 
any $f$ such that $f = \sum_{j=0}^\infty h_j$ holds in $\cs'(\rl)$ with
$\supp \cf_{n-1} h_j \subset \{ \xi' \mid |\xi'|_{a'} \le 2^j\}$.  
Analogously to $K$, the extension $Ef$ is defined
(whenever the following series converges in $\cs'$) as
\begin{equation}\label{3-1-10}
E f (x',x_n) = \sum_{j=0}^\infty \,  2^{-j a_n} \cfi_1 [\eta_j]
(x_n)\,  \, h_j (x') \, .
\end{equation}
Observe that this is not a mapping\,---\,despite the notation\,---\,since
each $f$ equals many sums like $\sum h_j$.

By the homogeneity properties of the Fourier transform,
\begin{equation}
\lim_{N \to \infty} \sum_{j=0}^N \, 2^{-j a_n} \cfi_1 [\eta_j]
(0)\,  \, h_j (x') 
= \sum_{j=0}^\infty \, h_j (x') = f(x') \, .  
\end{equation}
Hence $\g Ef=f$ is valid for such $f$, in particular when $f$ belongs to
some $A^{s,a'}_{p,q}(\Rl)$.
For simplicity's sake, a particular extension $Ef$
is said to depend {\em boundedly\/} on $f$ (although $Ef$ is not a map) when
\eqref{E'-eq} or \eqref{E''-eq} below holds:

\begin{satz}   \label{Ebd-thm}
There exists a constant $c$ such that the inequalities 
\begin{align}
  \|\, Ef \, | \bspq (\Rn)\| &\le c \, \| \, f \, |A^{s -
  \frac{a_n}{p},a'}_{p,q}(\Rl)\| ,
\label{E'-eq} \\
  \|\, Eg \, | \fspq (\Rn)\| &\le c \, \| \, g \, |A^{s -
  \frac{a_n}{p},a'}_{p,p}(\Rl)\| 
  \label{E''-eq}
\end{align}
hold for all $f \in A^{s - \frac{a_n}{p},a'}_{p,q}(\Rl)$
and all $g \in A^{s - \frac{a_n}{p}, a'}_{p,p}(\Rl)$, respectively.
\end{satz}

\begin{rem} Approximation spaces like $\aspq$ have been used for the trace
problem earlier, e.g. directly by P.~Oswald~\cite{Os}  and
in the proofs of V.~I.~Burenkov and M.~L.~Gol'dman \cite{BG}.
Traces of the approximation spaces themselves have been investigated by Yu.~V.~Netrusov~\cite{Ne}.
\end{rem}

\subsection{Consequences for the trace}
  \label{collect-ssect}

Using the tools from the previous subsection, one can now prove

\begin{satz}
  \label{aspq-thm}
The assertions of the introduction's main theorem are valid, and for any
$v(x')$ in $A^{s-\anp,a'}_{p,q}(\Rl)$ or $A^{s-\anp,a'}_{p,p}(\Rl)$ there is
an extension $Ev$ in $\bspq(\Rn)$ or $\fspq(\Rn)$, respectively, depending
boundedly on such $v$; in the generic cases $Ev$ may be defined by means of a
bounded operator $Kv$.
\end{satz}

It may be beneficial to list the consequences for the
various borderline cases, so we will formulate separate results with
detailed references. For brevity it is in the following understood that
the extensions $Ef$ depend boundedly on $f$, when $f$ is viewed as an element of
the pertinent range space for $\g$.

\begin{cor}   \label{Lp-prop}
Let $1 \le p < \infty$ and $0 < q \le 1$.  Then the trace operator
$\gamma_0$ maps $B^{\frac{a_n}{p},a}_{p,q}(\Rn)$ 
continuously onto $L_{p}(\rl)$ with extensions  
$Ef \in B^{\frac{a_n}{p},a}_{p,q}(\Rn)$
depending boundedly on $f$.
\end{cor}

The first contribution to this was made by S.~Agmon and L.~H\"ormander \cite{AgH},
who dealt with $p=2$, $q=1$ and $a=(1, \ldots,1)$.
For $1 \le p <\infty$,  $q=1$ and $a=(1, \ldots , 1)$, the first part of the
the result was stated by J.~Peetre \cite{Pe1}. The anisotropic variant 
was proved by V.~I.~Burenkov and M.~L.~Gol'dman \cite{BG} for $q=1$. Later 
M.~Frazier and B.~Jawerth \cite{FJ1} and J.~Johnsen \cite{Jo3} gave proofs
(except for the extensions) by other
methods for all $q\le 1$ in the isotropic situation. 

\begin{cor}   \label{Cb-prop}
When $0 < q \le 1$, then the trace $\gamma_0$ maps $B^{0,a}_{\infty,q}(\Rn)$
continuously onto $C_{\op{b}}(\rl)$ with extensions $Ef \in B^{0,a}_{\infty,q}(\Rn)$ 
depending boundedly on $f$.
\end{cor}

For $q=1$ the first proof of this result was given by V.~I.~Burenkov and
M.~L.~Gol'dman \cite{BG}.

\begin{cor}    \label{B-prop}
When $p$, $q\in\,]0,1]$, then $\gamma_0$ maps 
$B^{\frac{\a}{p} - |a'|,a}_{p,q}$
continuously onto the space ${A}^{|a'|(\ep - 1),a'}_{p,q}(\Rl)$ with extensions 
$Ef \in  B^{\frac{|a|}{p} - |a'|,a}_{p,q}(\Rn)$
depending boundedly on $f$.   
\end{cor}

Corollary~\ref{B-prop} should be a novelty; the determination of the trace
space as an approximation space does not, to our knowledge, have any forerunners.

\begin{cor}   \label{F1-prop}
The trace operator $\gamma_0$ maps $F^{a_n,a}_{1,q}(\Rn)$ continuously onto
$L_{1}(\rl)$, and there are extensions $Ef \in
F^{a_n,a}_{1,q}(\Rn)$ depending boundedly on $f$.
\end{cor}

In \cite{Ga}, Gagliardo proved $\gamma_0 (W^1_1 (\rn) )
= L_1 (\rl)$, which is closely related to Corollary~\ref{F1-prop} because of the embedding
$F^{1}_{1,2}(\rn) \hookrightarrow W^1_1(\rn) \hookrightarrow B^1_{1,\infty}(\rn)$.  
M.~Frazier and B.~Jawerth \cite{FJ2} were the first to attempt a proof of the
corollary, but their argument seems rather flawed;
cf.\ Remark~\ref{FJ-rem} below; H.~Triebel obtained the first part of Corollary~\ref{F1-prop}
by another approach based on atomic decompositions \cite[4.4.3]{Tr9}, but without making it
clear that $F^{a_n,a}_{1,\infty}(\Rn)$ is a subspace of $C(\re,\cd'(\rl))$.

The next result is an analogue of Corollary~\ref{B-prop}.

\begin{cor}   \label{F-prop}
When $0 < p < 1$,  then $\gamma_0$ maps 
$F^{\frac{|a|}{p} - |a'|, a}_{p,q}(\Rn)$
continuously onto the space ${A}^{|a'|(\ep - 1),a'}_{p,p}(\Rl)$ with
extensions $Ef \in F^{\frac{|a|}{p} - |a'|, a}_{p,q}(\Rn)$ depending
boundedly on~$f$.           
\end{cor}

A final remark to the existence of a linear right inverse of $\gamma_0$:
J.~Peetre \cite{Pe1} ($a=(1,\ldots, 1)$) and V.~I.~Burenkov and M.~L.~Gol'dman \cite{BG}
(general anisotropic case) have shown that if $1 \le p < \infty $ and
$s=a_n/p$, then there exists no linear extension operator mapping
$L_p (\rl)$ boundedly  to $B^{\frac{a_n}{p},a}_{p,1}(\Rn)$;
whether this remains true in other borderline cases seems to be unknown.

\begin{rem}
  \label{FJ-rem}
Directly below \cite[Th.~11.1]{FJ2} the authors write (in our notation): ``We will show
directly that $\g(F^{s}_{p,q})$ is independent of $q$. Given this, we have
$\g(F^s_{p,q})=\g(F^s_{p,p})$ and all conclusions follow from the [Besov space] results in
\cite[Sect.~5]{FJ1}, since  $B^{s}_{p,p}=F^{s}_{p,p}$.''.  However, it is
evident from their proof that they tacitly assume $\g$ to be defined on $F^s_{p,q}$ for
{\em two\/} arbitrary sum-exponents ($q$ and $r$ in $\,]0,\infty]$), and not
just for $q=p$, but they never
support this implicit claim by arguments.

Although it is true (and trivial to verify for the generic cases as
well as for $s=\anp-\aa$ when $p<1$, in view of the Sobolev embedding into
$B^{a_n,a}_{1,1}(\Rn)$), it does require a proof that $\g$ is well defined on
$F^{a_n,a}_{1,q}(\Rn)$, since if $1<q\le\infty$ this space is strictly
larger than any Besov space on which $\g$ is defined; cf.\ the vertex in
Figure~\ref{param-fig}. 
The authors claim to have covered such $F$-spaces
as a novelty, but in view of the described flaw it should be appropriate
that we  prove 
that $F^{a_n,a}_{1,\infty}(\Rn)\hookrightarrow C_{\op{b}}(\re,L_1(\rl))$ and
that $\g (F^{a_n,a}_{1,q})$ is independent of $q\in\,]0,\infty]$;
see Proposition~\ref{Cembd-prop} ff.\ and Proposition~\ref{indep-prop} below.
\end{rem}



\section{A sharper comparison of the Besov--Lizorkin--Triebel classes and
  the approximation spaces}  \label{apr-sect}

In Proposition~\ref{aspq-prop} we have identified $\aspq$ with standard
function spaces 
in the generic cases. We now investigate the remaining borderline case 
$s= |a|(\ep -1)$ for $p<1$; the analysis involves the
continuity properties of $\g$ proved in Theorem~\ref{aspq-thm} above. 

Concerning the borderline with $s=\a(\ep-1)$, it is noteworthy that the two cases
$q \le p <1$ and $p < q \le 1$ give quite different results:

\begin{satz}   \label{apr-thm}
Let $0 < p < 1 $ and $0 \le q \le 1$. 
Then
\begin{equation}\label{5-2-1}
A^{\a(\ep - 1),a}_{p,q} (\Rn) \hookrightarrow  
B^{\a(\frac 1 r  -1),a}_{r,u}(\Rn) 
\end{equation}
holds if, and only if, 
\begin{equation}\label{5-2-2}
\max (p,q) \le r \le \infty \qquad \mbox{and} \qquad u=\infty \, ,
\end{equation}
and whenever $0 < r < \infty$ and $0 < u \le \infty$, then
\begin{equation}\label{5-2-3}
A^{\a(\ep - 1),a}_{p,q} (\Rn) \not\subset  
F^{\a(\frac 1 r  -1),a}_{r,u}(\Rn).
\end{equation}

Conversely, 
\begin{equation}\label{5-2-5}
B^{\a(\frac 1 r - 1),a}_{r,u}(\Rn) \hookrightarrow 
A^{\a(\ep - 1),a}_{p,q} (\Rn)
\end{equation}
holds if
\begin{equation}\label{5-2-6}
0 < r \le p \qquad \mbox{and} \qquad 0 < u \le q \, , 
\end{equation}
while
\begin{equation}\label{5-2-7}
F^{\a(\frac 1 r - 1),a}_{r,u} (\Rn) \hookrightarrow
A^{\a(\ep - 1),a}_{p,q}(\Rn)
\end{equation}
holds if one of the following conditions does so:
\begin{equation}\label{5-2-8}
 \bigl( 0 < r < p \quad\text{and}\quad  r\le q \bigr) 
\qquad \mbox{or} \qquad \big( r = p\le q \quad \mbox{and} \quad
0 < u \le q \big).
\end{equation}
\end{satz}

The necessity of the conditions \eqref{5-2-6} and \eqref{5-2-8} has been
obtained for the parts concerning $r$ and $p$, but not for the sum-exponents;
cf. Remark~\ref{ru-rem} below.

By application of the above results, it is clear that, for $p < 1$,
one has embeddings
\begin{equation}
    \label{bab-eq}
B^{\a(\ep - 1),a}_{p,q}(\Rn ) \hookrightarrow 
A^{\a(\ep - 1),a}_{p,q}(\Rn) \hookrightarrow 
B^{\a(\frac1r -1),a}_{r,\infty}(\Rn) , \qquad r=\max(p,q).
\end{equation}
The latter is \emph{optimal} in the sense that any space $\bspq$ or
$\fspq$ (with $s=\a(\ep-1)$), which $\aspq$ is embedded into, 
actually also has the Besov space on the
right hand side as an embedded subspace. This follows from
\eqref{5-2-2}--\eqref{5-2-3} and the usual embeddings.

However, a somewhat sharper argument yields the following result:

\begin{prop}
  \label{ABF-prop}
For $0<p<1$ and $0<q\le1$ the space $A^{\a(\ep-1),a}_{p,q}(\Rn)$ is neither
a Besov space $B^{s',a}_{p',q'}(\Rn)$ nor a Lizorkin--Triebel space
$F^{s',a}_{p',q'}(\Rn)$ for any admissible $(s',p',q')$. 
\end{prop}

From this proposition, from Theorem~\ref{aspq-thm}, and
the fact that also $L_1$ is neither
a Besov nor a Lizorkin--Triebel space, we obtain 

\begin{cor}
For $0 < p \le 1$ the ranges of the trace operator $\g$,
\[\gamma_0 \, (B^{\frac{|a|}{p}-|a'|,a}_{p,q}(\Rn))\, , 
\quad 0 < q \le 1\, , \qquad \mbox{and} \qquad
\gamma_0 \, (F^{\frac{|a|}{p} - |a'|,a}_{p,q}(\Rn))\, , 
\quad  0 < q \le \infty\, , 
\] 
do not belong to the scales of anisotropic 
Besov--Lizorkin--Triebel spaces on $\re^{n-1}$.
\end{cor}

The embedding of the trace spaces 
$\gamma_0 (B^{\frac{|a|}{p}-|a'|,a}_{p,q}(\Rn))$ into 
$B^{|a'|(\frac{1}{r}-1),a'}_{r,\infty}(\Rl))$ for $r= \max (p,q)$,
cf. (\ref{5-2-1}), was first proved by Johnsen \cite{Jo3} 
(isotropic situation). By the above discussion, the present results are
sharper and optimal.

A final remark on the borderline cases. For fixed $p$ one may ask for the
largest spaces on which the trace exists. One should clearly minimise 
$s$ and maximise $q$ (which is done throughout this paper), but the
dependence on the anisotropy $a$ may also be considered. For the isotropic
spaces (indicated without the $a$) the following embeddings hold, for $p<\infty$:
\begin{gather}\label{ae}
  B^{\ep + (n-1)(\ep -1)_+}_{p,q}(\rn) \hookrightarrow 
  B^{\frac{a_n}{p} + |a'|(\ep -1)_+,a}_{p,q} (\Rn),
\\
  F^{\ep + (n-1)(\ep -1)_+}_{p,q}(\rn) \hookrightarrow 
  F^{\frac{a_n}{p} + |a'|(\ep -1)_+,a}_{p,q} (\Rn);
\end{gather}
i.e.\ the anisotropic spaces are larger than the corresponding isotropic ones. 
Furthermore, if $p < 1$, then the Besov spaces  
$B^{\frac{a_n}{p} + |a'|(\ep -1)_+,a}_{p,1} (\Rn)$ and  
the Lizorkin--Triebel spaces $F^{\frac{a_n}{p} + |a'|(\ep -1)_+,a}_{p,\infty}(\Rn)$ 
are the largest possible, however, all of them are contained in
$F^{a_n ,a}_{1,\infty} (\Rn)$ which is the largest possible for $p=1$. 
If $1 < p < \infty$, the only candidate is the Besov space
$B^{\frac{a_n}{p},a}_{p,1} (\Rn)$.   
Returning to the anisotropic classes with different 
anisotropies, it is not hard to see that the spaces 
depend increasingly on each component $a_j$ of $a$. 
So, although  $a_n$ is bounded by $a_n \le s \cdot p$, 
within the anisotropic spaces there is no maximal spaces on which the trace makes
sense.
For $p=\infty$ there is a largest space having a continuous 
trace, namely $C_{\op{b}}(\rn)$.
But, as mentioned in the introduction, all these anisotropic spaces
are subclasses of $C(\re, \cd'(\rl))$ with its natural notion of a trace.


\section{Remaining Proofs}
  \label{proof-sect}
Some proofs below are essentially just anisotropic variants of known
techniques (scattered in the journals and main references like
\cite{St,Tr6, Tr9}), but even so the most important ones are presented (or
sketched) 
here for the reader's convenience. It would lead too far to do this
consistently, so in the remaining cases we shall have to make do with indications
of the necessary changes, however. 


\subsection{The assertions in Section~\ref{gen-sect}} 
\subsubsection{Proof of Theorem~\ref{gen-thm}} 
  \label{gen-ssect}
By means of elementary embeddings (cf.\ around \eqref{700-8} below),
the `only if' part is, for $p<\infty$, a consequence of the 
following lemma, which is carried over from \cite[Lem.~2.8]{Jo2}.

\begin{lem}  \label{711-l4}
Let $0<p<\infty$. For any $u\in \cs (\rl )$ there exists 
$u_k\in \cs(\rn )$
such that
\begin{gather}
  \g u_k (x')=u(x') \quad \mbox{for any}\quad k\in \N , 
\\
  \lim\limits _{k\to \infty} u_k=0\quad \mbox{in}\quad 
       B_{p,q}^{\frac{a_n}{p},a}(\rn )\quad \mbox{if}\quad 1<q\le \infty  ,
  \label{715-4} 
\\
  \lim\limits _{k\to \infty} u_k=0\quad \mbox{in}\quad 
       F_{p,q}^{\frac{a_n}{p},a}(\rn )\quad \mbox{if}\quad 1<p< \infty  .
  \label{715-5}
\end{gather}
If $0<p\le1$, there exists $v_k\in \cs(\rn )$ such that 
$\lim_{k\to \infty}  \g v_k=\delta _{0}$ in $\cs'(\rl)$ while
\begin{equation}
  \label{715-7}
  \lim\limits _{k\to \infty} v_k=0 \quad \mbox{in}\quad 
  B_{p,q}^{\frac{|a|}{p}-|a'|,a}(\rn )
  \quad\mbox{for}\quad 1<q\le \infty ;
\end{equation}
here $\delta_0$ stands for the Dirac measure at $x'=0$.
\end{lem}
To prove Lemma \ref{711-l4} one can set $u_k(x)=u(x')w_k(x_n)$ with
\begin{equation}
  w_k (x_n)=\frac{1}{k}\sum\limits _{l=1}^k w(2^{la_n} x_n)  
\end{equation}
where $w\in \cs (\re )$ with $\supp \widehat{w} \subset \{\,\xi _n\in \re \mid
\frac{3}{4}\le |\xi _n|^{1/a_n}\le 1\,\}$ and $w(0)=1$.
Moreover,
\begin{equation}
   v_k(x',x_n)=\frac{1}{k}\sum\limits _{l=k+1}^{2k}
   2^{l|a'|}\, f(2^{la'}x')\, g(2^{la_n}x_n)
\end{equation}
has the claimed properties at least when $f\in \cs (\rl )$ and $g\in \cs (\re )$ satisfy
\begin{gather}
  \supp \cf _{n-1}f\subset \bigl\{\,\xi '\in\rl \bigm|
   |\xi '|_{a'}\le \tfrac{1}{2}\,\bigr\} \quad \mbox{and}\quad 
  \int _{\rl} f(x')\, dx' =1,
\\
   \supp \cf g\subset \bigl \{\,
            \xi _n\in \re \bigm|
   |\xi _n|^{1/a_n}\le \tfrac{1}{2}
                       \,\bigr\}\quad \mbox{and}\quad
   g(0)=1.
\end{gather}
Indeed, it is not hard to check that the norms are ${\mathcal
O}(k^{\tfrac{1}{q}-1})$, respectively ${\mathcal O}(k^{\ep-1+\varepsilon})$ for
the $\fspq$-norm,
see \cite[Lem.~2.8]{Jo2}, for one may use the fact that if
$s>0$ (it is here the restriction $p<\infty$ 
is needed), then there exists a constant $c>0$ such that 
\begin{gather}
  \label{715-1}
  \|f\otimes g\, |\, \bspqa (\rn )\|\le c\, 
  \|f\, |\, B_{p,q}^{s,a'}(\Rl )\|\, \|g\, |\, B_{p,q}^{s/a_n}(\re )\| ,
\\
   \label{715-01}
  \|f\otimes g\, |\, \fspqa (\rn )\|\le c\, 
  \|f\, |\, F_{p,q}^{s,a'}(\Rl )\|\, \|g\, |\, F_{p,q}^{s/a_n}(\re )\|.
\end{gather}
The isotropic version of \eqref{715-1}--\eqref{715-01} is due to J. Franke, see
\cite[Lem.~1]{Fr}. For the anisotropic version one may use a
paramultiplicative decomposition (in Yamazaki's sense) of the direct product
$f\otimes g$ and apply the estimates in \cite{Ya} (cf.\ the
case $f=\delta_0$ treated in \cite{Jo2}).

For the space $B^{0,a}_{\infty,q}(\Rn)$ with $1<q\le\infty$, one can modify the
proof of the lemma by taking $u\in\cs$ to have a sufficiently small (non-empty) spectrum
such as the ball $\{\,\,\xi'\mid |\xi'|_{a'}\le \tfrac{1}{2}\,\}$. Then the
norm of $u_k$ is seen to be ${\mathcal O}(k^{\eq-1})$ if, instead of
\eqref{715-1}, one calculates directly by means of an anisotropic Lizorkin
representation (in the language of \cite[2.6]{Tr6}) with a \emph{smooth} partition
of unity. (This means that the $\ffij$ entering the norm of $\bspq$
should be replaced by some $\eta_{jk}$, with $k$ in a finite $j$-independent
set, such that $\cup_k \supp\eta_{jk}$ equals a `corridor' (the complement
of an $n$-dimensional rectangle inside a dilation by $2^{a}$ of itself) and
such that each $\eta_{jk}$ is a product
$\theta_{k}(2^{ja'}\xi')\omega_{k}(2^{ja_n}\xi_n)$; this is well known to
give an equivalent quasi-norm for $\bspq$ by Lemma~\ref{cnv1-lem}; an isotropic
version may be found in [Jo3].)

\bigskip

It remains to prove the stated continuity of
$\g$. For the Lizorkin--Triebel case one may show that the operator $\g$ is
defined on $\fspq$ for every $(s,p,q)$ considered in the theorem, see
Proposition~\ref{Cembd-prop}, and that $\g(\fspq)$ is independent of $q$; this last fact is
obtained in Appendix~\ref{qindep-app} below. Then the identification 
$B^{s,a}_{p,p}=F^{s,a}_{p,p}$ reduces the question to the Besov case; for
the boundedness of $\g$ one may use the inequality \eqref{qt-est} in
Appendix~\ref{qindep-app} below. 

For the treatment of $\bspq$, we use the short argument of
\cite[Sect.~3]{Jo3}; the idea is to combine the Nikol'skij--Plancherel--Polya
inequality with the Paley--Wiener--Schwartz Theorem to deduce the crucial
mixed-norm estimate \eqref{basic-ineq}.

First we  remark that if $h \in \cs ' (\rn )$ satisfies 
$\supp \cf h \subset \{\, \xi\in \rn\mid |\xi|_a \le A \,\} $ 
for some $A>0$, then the restriction  
$h(x',\cdot)$, obtained by freezing $x'$, fulfils 
\begin{equation}
  \label{713-21}
  \supp  \cf_{x_n\mapsto \xi _n} h(x',\xi _n) \subset
  \{\, \xi _n\in \re \mid | \xi _n| \le A^{a_n} \,\} .
\end{equation}
This is a consequence of \eqref{112-6} below and the 
Paley--Wiener--Schwartz Theorem (cf.~\cite[7.3.1]{hoer}), for these give
that $g:=h(x', \cdot )$ is analytic and satisfies  
\begin{equation}
  g(x_n+iy_n)\le C(x')(1+|x_n|)^{N}e^{A^{a_n}|y_n|}. 
\end{equation}

Now, let $f\in\bspq(\Rn)$. Applying (\ref{713-21}) with $A=2^{j}$ to 
$\cf ^{-1}[\ffij \cf f](x', \cdot )$, the
Nikol'skij--Plancherel--Polya inequality, see for example
\cite[1.3.2]{Tr6}, yields
\begin{equation}
  \label{713-22}
  \sup _{x_n \in \re} | \cf ^{-1} [\ffij\cf f](x',x_n) |
  \le c \, 2^{ja_n/p} \big (\int_{\re}  
  | \cf ^{-1} [\ffij \cf f](x',x_n) |^p \, dx_n  \big )^{1/p} 
\end{equation}
where the constant $c$ does not depend on $x'$, $f$ and $j$;
hence $x'$-integration gives
\begin{equation}
  \label{basic-ineq}
  \| \sup_{x_n} \cf ^{-1} [\ffij\cf f] (\cdot,x_n)\, |\,  L_p(\rl )\|
 \le 
  c \, 2^{ja_n/p}\, \| \cf ^{-1} [\ffij \cf f] \,|\ L_p(\rn )\|  .
\end{equation}
Estimating the supremum from below by the value for $x_n=0$ and 
arguing as for \eqref{713-21},  we obtain
\begin{equation}
  \label{basic-eq}
  \supp  \cf_{x'\mapsto \xi'}\bigl( \cf ^{-1} [\ffij\cf f](x',0)\bigr) \subset
  \{\, \xi'\in \rl \mid | \xi'|_{a'} \le 2^{j} \,\} .
\end{equation}
So under the restriction $s>\anp+\aa(\ep-1)_+$, Lemma~\ref{cnv2-lem} yields
the convergence of the series $\sum_{j=0}^\infty \cf ^{-1} [\ffij \cf f] (\cdot, 0)$
in $\cs '(\rl )$ as well as the boundedness of $\gamma_0 \colon \bspq(\Rn )
\to B^{s-\frac{a_n}{p},a'}_{p,q}(\Rl)$.

\subsubsection{Proof of Proposition~\ref{Cembd-prop}}
  \label{Cembd-pf}
The particular value $x_n=0$ is unimportant in \eqref{basic-eq}, so when
the application of the Nikol'skij--Plancherel--Polya inequality above 
is repeated, a slightly stronger conclusion is reached: for $p_1:=\max(1,p)$
and $f_j:=\cfi[\ffij\cf f]$,
\begin{equation}
   \sup_{x_n\in\re} \|  f_j (\cdot, x_n)\, |\,  L_{p_1}(\rl )\|
 \le 
  c 2^{j(\tfrac{a_n}{p}+\aa(\ep-1)_+)}\, \|\, f_j\,
  | L_p\|  .
\end{equation}
If $C_{\op{b}}$ temporarily stands for bounded, continuous functions, the
left hand side is the norm of $f_j$ in $C_{\op{b}}(\re,L_{p_1}(\rl))$, and 
this is in $\ell_{1}(\n)$ with respect to $j$ because the right hand side is
so. Consequently, the series $\sum f_j$ converges in
$\cs'(\rn)$ to $f$ and in  $C_{\op{b}}(\re,L_{p_1}(\rl))$, and because the
latter space is continuously embedded into the former (a reference to this
folklore could be Prop.~3.5 and (5.4) in \cite{Jo3}), this shows that $f\in
C_{\op{b}}(\re,L_{p_1}(\rl))$. Moreover, the triangle inequality applied to
$f=\sum f_j$ yields continuity of 
\begin{equation}
\begin{split}
 \bspq(\Rn)\hookrightarrow  C_{\op{b}}(\re,L_{p_1}(\rl)) &\qquad\text{whenever
 $s\ge \tfrac{a_n}p+\aa(\tfrac1p-1)_+$,}
\\
 &\qquad\text{provided $q\le1$ in case of equality}.
\end{split}  
\end{equation}
Finally, the uniform continuity with respect to
$x_n$ should be verified. However, if translation by $h\in\re$ is denoted by
$\tau_h$, that is $\tau_hf(x)=f(x',x_n-h)$, then the boundedness above yields
the following, say with $s=\anp+\aa(\ep-1)_+$ for simplicity:
\begin{equation}
  \sup_{x_n\in\re}
  \|\,{f(\cdot,x_n)-f(\cdot,x_n-h)}\,|{L_{p_1}(\rl)}\|\le
  c(\sum_{j=0}^{\infty}2^{sjq}\|\,{(1-\tau_h)f_j}\,|{L_{p_1}}\|^q)^{\frac1q}.
\end{equation}
Indeed, this is clear since $1-\tau_h$ commutes with $\cfi\ffij\cf$, and so it
remains to note that the right hand side tends to $0$ for $h\to 0$ by
majorised convergence.

Since $\fspq\hookrightarrow B^{s,a}_{p,\infty}$ may be used for the
generic Lizorkin--Triebel cases, it suffices to consider $\fspq$ in the 
borderline cases with $s=\ap-\aa$ and $p\le1$; however, by the  Sobolev
embeddings it is enough to treat $p=1$, hence to show that 
\begin{equation}
    \label{F11-eq}
  F^{a_n,a}_{1,\infty}(\Rn)\hookrightarrow C_{\op{b}}(\re,L_1(\rl)).
\end{equation}
To do so, we replace the above use of the
Nikol'skij--Plancherel--Polya inequality by an application of the Jawerth
embedding $F^1_{1,\infty}(\re)\hookrightarrow B^0_{\infty,1}(\re)$. As above,
for $f\in F^{a_n,a}_{1,\infty}(\Rn)$,
\begin{equation}
  f(x)=\sum_{j=0}^\infty f_j(x).
  \label{LP-eq}
\end{equation}
Note, as a preparation, that this series converges pointwise for
$x\notin N$, where $N$ is a Borel set in $\rn$ with $\op{meas}(N)=0$; indeed,
this follows since $\sum \|\,{f_j}\,|{L_1}\|<\infty$ must hold for any $f$ in
$F^{a_n,a}_{1,\infty}(\Rn)$. 

By Fubini's theorem, there is also a null set
$M\subset\rl$ such that $x'\notin M$ implies
\begin{equation}
  \int_{\re} \sup_j |2^{ja_n}f_j(x',x_n)|\,dx_n<\infty.
  \label{Fint-eq}
\end{equation}
Invoking Lemma~\ref{cnv2-lem} one therefore obtains a function $x_n\mapsto
g(x',x_n)$ in $F^{a_n}_{1,\infty}(\re)$ for which $g(x',\cdot)=\sum
f_j(x',\cdot)$; using $a_n\ge1$ to apply the Jawerth embedding, we have
\begin{equation}
  \begin{split}
  |g(x',x_n)| &\le \|\,g(x',\cdot)\,| B^{0}_{\infty,1}(\re)\|
\\
  &\le c \|\,\sum f_j(x',\cdot)\,| F^{a_n}_{1,\infty}(\re)\|
   \le c'  \smash[t]{\int_{\re}} \sup_j |2^{ja_n}f_j(x',x_n)|\,dx_n<\infty.
  \end{split}
  \label{gf-ineq}
\end{equation}
The $x'$-dependence is not arbitrary as it seems to be, for there is another
null set $M'\subset\rl$ such that $M\subset M'$ and when $x'\notin M'$, then
\begin{equation}
  g(x',x_n)= f(x',x_n)\quad\text{for $x_n$ a.e.\ in $\re$}.
  \label{gf-eq}
\end{equation}
Indeed, the section $N_{x'}=\{\,x_n\mid (x',x_n)\in N\,\}$ is
a Borel set and the relation
$0=\op{meas}(N)=\int_{\rl}\op{meas}(N_{x'})\,dx'$ gives that
$\op{meas}(N_{x'})=0$ for $x'$ outside a null set $M'$, which may be
assumed to contain $M$. So by \eqref{LP-eq}, $f(x',\cdot)=\sum f_j(x',\cdot)$ holds
outside $N_{x'}$ whenever $x'\notin M'$. But, since the norm series $\sum
\|\,f_j(x',\cdot)\,| L_1(\re)\|$ is 
estimated by the integral in \eqref{Fint-eq}, hence is finite, the series
$\sum f_j(x',\cdot)$ converges to $g(x',\cdot)$ in $L_1(\re)$. So, by the
fact that a pointwise limit coincides a.e.\ with a limit in mean,
\eqref{gf-eq} is obtained.

It is thus justified to integrate both sides of \eqref{gf-ineq} with respect
to $x'$, and therefore 
\begin{equation}
  \sup_{x_n\in\re}
  \|\,{f(\cdot,x_n)}\,|{L_1(\rl)}\|\le
  c\|\,{f}\,|{F^{a_n,a}_{1,\infty}}\|.
\end{equation}
Now the uniform continuity in $x_n$ of any $f$ in $F^{a_n,a}_{1,\infty}$ may
be shown by an argument analogous to the Besov case above. This proves
\eqref{F11-eq} and thus the proposition.

\subsubsection{Proof of Theorem~\ref{gen1-thm}}   \label{gen1-sssect}
For part (i), the streamlined method of \cite[Sect.~3]{Jo3} may be adopted
as follows. Writing
\begin{equation}
K v(x',x_n)= \smash{\jsum} u_j(x',x_n)  ,
\end{equation}
where for any $j\in \no$
\begin{equation}
  u_j(x',x_n)= 2^{-ja_n} \, (\cf _{1}^{-1}\eta _{j})(x_n)\,
  \cf _{n-1}^{-1}[\ffij (\cdot,0)\, \cf _{n-1}v](x'),  
\end{equation}
it is straightforward to see that $u_j$ has a compact spectrum:
using the triangle inequality one can see that 
$|\xi |_a\sim |\xi '|_{a'}+|\xi_n|^{1/a_n}$ and then find a constant 
$A>0$ such that
\begin{equation}
    \supp \cf u_j \subset\{\,\xi \mid 2^{ja_n}\le 
    |\xi _n|\le 2^{(j+1)a_n},\ |\xi '|_{a'}\le c\, 2^j \,\}
    \subset 
    \{\,\xi \mid \tfrac{2^j}{A}\le|\xi |_{a}\le A\, 2^j\, \}.
  \end{equation}
For any $j\in \no $ we have, when $\check\eta:=\cf _{1}^{-1}\eta$ and $\hat v:=\cf _{n-1}v$, 
\begin{equation}
  2^{sj}\|\,u_j\, | L_p(\rn)\| = 2^{j(s-\anp)} 
      \|\,\check\eta \, | L_p(\re )\|  
  \cdot 
 \|\,\cf _{n-1}^{-1}[\ffij (\cdot, 0)\hat v] \,| L_p(\rl )\|  
\end{equation}
where the right hand side is in $\ell_q$ provided $v\in B^{s-\anp,a'}_{p,q}(\Rl)$.
By Lemma \ref{cnv1-lem}, this implies 
\begin{equation}
  \|\, K v\, |\bspqa (\rn )\| \le 
  c \bigl ( \jsum
         2^{jsq}\|\, u_j\,| L_p(\rn )\|^{q}
   \bigr)^{\frac1q} 
  \le c' \|\,v\, | B_{p,q}^{s-\frac{a_n}{p},a'}\|
\end{equation}
and the proof of (i) is complete.
 
\bigskip

When $v\in B^{s,a'}_{p,p}(\Rl)$ one can modify the corresponding proof of
\cite[2.7.2]{Tr6} in the following way: by the $\fspq$-part of
Lemma~\ref{cnv1-lem}, the series for $Kv$ converges to an element of
$F^{s+\anp,a}_{p,q}(\Rn)$ if we can show a certain estimate; this is done as
in \cite{Tr6}. The main thing is to get the correct auxiliary inequality
which is
\begin{equation}
  |2^{-ja_n}\cfi_{1}\eta_j(x_n) | \le c(1+2^{ja_n}|x_n|)^{-\delta}
\end{equation}
for a sufficiently large positive $\delta$; in addition it is convenient to
split the integration there over the subintervals $I_l=\,]-2^{-la_n},-2^{-(l+1)a_n}]\cup
[2^{-(l+1)a_n},2^{-la_n}[$.


\subsection{Assertions in Section~\ref{bdr-sect}}
  \label{bdr-ssect} 
\subsubsection{Proof of Proposition~\ref{aspq-prop}}

For part (i) the definitions and line \eqref{700-0} below yield 
$\bspq (\Rn) \hookrightarrow \aspq (\Rn)$.
However, if $s>\a(\ep - 1)_+$, then Lemma~\ref{cnv2-lem}, cf. the appendix,
implies that the two spaces are equal. 

To prove parts (ii) and (iii), let $f \in L_p (\rn )$ for some $p\in
[1,\infty[\,$ (the case $f\in C_{\op{b}}(\rn)$ for $p=\infty$ is treated similarly).
Selecting a subsequence $j_k$ of $\no$ such that $\lambda_k:=2^{j_k}$ (with
$\lambda_0=1$) satisfies
\begin{equation}
\|\, \cf ^{-1} [\psi (\lambda_k^a\cdot)\cf f] - f\, | L_p (\rn )\| < 
2^{-k-1} \, \| \, f \, | L_p(\rn)\|
\qquad\text{for }  k\in \n
\end{equation}
(which is possible as seen from the usual convolution estimates), one
clearly has that the series
\begin{equation}
  f(x) = \cf ^{-1} [\psi\cf f](x) + 
  \sum_{k=1}^\infty \, \big(\cf ^{-1} [\psi (\lambda_{k}^{a}\cdot)\cf f](x)
  -\cfi [\psi (\lambda_{k-1}^a\cdot)\cf f](x)\big)
  \label{fapp-eq}  
\end{equation}
converges in $\| \, \cdot \, | L_p (\rn )\|$. Setting $h_j$ equal to the
$k^{\op{th}}$ summand in \eqref{fapp-eq} for $j=j_k$ and $h_j=0$ for all
other $j$, it is found by the definition of the $\aspq$ that $f \in
A^{0,a}_{p,q}(\Rn)$ for all $q$ and  
\begin{equation}
\| \, f \, | A^{0,a}_{p,q}(\Rn)\| \le  c \, \| \, f \, | L_p(\rn )\|.  
\end{equation}
The converse inclusion may be shown for $0<q\le1$, for when
$f\in A^{0,a}_{p,q}(\Rn)$ is written as $f=\sum h_j$ according to the definition, then
the embedding $l_q\imb l_1$ yields
\begin{equation}
  \jsum \|h_j\, |\, L_p(\rn )\|\le (\jsum
  2^{j\cdot 0\cdot q} \|h_j\, | L_p(\rn )\|^q)^{1/q}  ,
\end{equation}
so that the  completeness of $L_p$ gives
$\|f\, | L_{p}(\rn )\|\le \|f\, | A^{0,a}_{p,q}(\Rn )\|$.
The proof is complete.

\subsubsection{Proof of Theorem~\ref{Ebd-thm}}

Here one can use the same strategy as for Theorem~\ref{gen1-thm}, except
that a given $v$ in the approximation space should be written as $v=\sum h_j$
and $h_j$ should then replace $\cfi_{n-1} \ffij(\cdot,0)\cf_{n-1}
v$; this works because also $h_j$ has its spectrum in the ball
$\{\,\xi'\mid |\xi'|_{a'}\le 2^{j}\,\}$ and because the representation $\sum
h_j$ can be chosen such that its relevant norm is less than $2\|\,v\,|
A^{s-\anp,a'}_{p,q}(\Rl)\|$. 

\subsubsection{Proof of Theorem~\ref{aspq-thm} (and of the main theorem)}

Obviously (b) entails (a) in the main theorem; cf.\
Corollary~\ref{aspq-cor}. The fact that (a) implies (c) is proved in
connection with Theorem~\ref{gen-thm}, and the extensions depending
boundedly on $v$ were established in Theorem~\ref{Ebd-thm} and, for $K$, in
Theorem~\ref{gen1-thm}. So it remains to prove (c)$\implies$(b).

When $f$ is in $\bspq$ one derives \eqref{basic-ineq} as before;
since $s\ge \anp+\aa(\ep-1)_+$ it is straightforward to see from
\eqref{basic-ineq} that $\sum \cfi[\ffij \cf f](\cdot,0)$ converges in $L_p$
if $p\ge1$ or, by the Nikol'skij--Plancherel--Polya inequality, in $L_1$
if $p<1$. So we may write $\g f=\sum \cfi[\ffij \cf f](\cdot,0)$, and
\eqref{basic-ineq}, \eqref{basic-eq} and the definition of the approximation
spaces thereafter show the boundedness of $\g$ from $\bspq$ into
$A^{s-\anp,a'}_{p,q}$. The Lizorkin--Triebel case follows from the Besov
case as before; cf.\ Propositions~\ref{Cembd-prop} and \ref{indep-prop} below.

\subsection{The assertions in Section 4}
\subsubsection*{Proof of Theorem~\ref{apr-thm}}
Step~1. 
To deduce all the embeddings, note that when $s=\a(\ep-1)$ and
$r=\max(p,q)$, then it is easy to establish that 
\begin{equation}
  \bspq(\Rn)\hookrightarrow \aspq(\Rn) \hookrightarrow B^{\a(\frac1r-1),a}_{r,\infty}(\Rn);
\end{equation}
in fact, the first inclusion is obvious from the definitions, and the second follows
if Lemma~\ref{cnv3-lem} below is invoked in addition. Then the sufficiency
of \eqref{5-2-2} and \eqref{5-2-6} is clear in view of the Sobolev
embeddings, and from \eqref{700-8} or the anisotropic Jawerth embedding in
Proposition~\ref{711-t21} below it is seen analogously that \eqref{5-2-8}
implies \eqref{5-2-7}.  

Step~2. 
To show that \eqref{5-2-1} implies the $u$-part of \eqref{5-2-2} in
Theorem~\ref{apr-thm}, we shall use the already proved boundedness of $\g$;
it is therefore convenient to replace \eqref{5-2-1} by the assumption that
$A^{\aa(\ep-1),a'}_{p,q}(\Rl)\imb B^{\aa(\frac 1r-1),a'}_{r,u}(\Rl)$ for
some $n\ge 2$.

It suffices to find Schwartz functions $g_k$ such that
\begin{gather}
\label{717-1}
 \sup_{k\in \N} \|\,g_k\, | B_{p,q}^{\frac{|a|}{p}-|a'|,a}(\rn )\| <\infty,
\\
 \label{717-2}
 \|\,\g g_k\, |  B_{ru}^{|a'|(\frac{1}{r}-1),a'}(\rl )\|\ge
 c k^{1/u}\quad \text{for any}\quad k\ge 2 .
\end{gather}
Indeed, it then follows from \eqref{5-2-1} and the boundedness of $\g$ from
$B^{\ap-\aa,a}_{p,q}$ that
\begin{eqnarray}
 \label{717-3}
 ck^{1/u}                         
  \le c' \|\,\g g_k\, |  A^{|a'|(\frac{1}{p}-1),a'}_{p,q}(\rl )\|
   \le c''  \|\,g_k\, | B^{\frac{|a|}{p}-|a'|,a}_{p,q}(\rn )\|   ,
\end{eqnarray}
which contradicts \eqref{717-1} unless $u=\infty$.

To prove the existence of such $g_k$, note first that it is possible to take
the partition of unity $1=\sum \ffij$ such that, say, $\ffi_0(\xi)=1$ for
$|\xi|_{a}\le \frac{11}{10}$ and $\ffi_0(\xi)=0$ for $|\xi|_a\ge
\frac{13}{10}$ (this choice is consistent with the conventions in
\cite{Jo1}, which will be convenient later). Similarly the $\eta$ of
Section~\ref{K-ssect} should fulfil $\supp\eta\subset\,]1,\frac{21}{20}[$. 

With some $k_0$ to be determined, we set
\begin{equation}
\label{717-5}
  g_k(x',x_n)= (\cf _1^{-1}\eta)(2^{(k+k_0)a_n}x_n) \sum _{j=0}^{k} 
                \cf _{n-1}^{-1} [\ffij (\cdot,0)](x').
\end{equation}
Then $\supp\cf g_k$ is contained in the set where both
$|\xi '|_{a'}<\frac{13}{10}\cdot 2^{k}$ and 
$2^{(k+k_0)a_n}< |\xi_n|<\frac{21}{20}\cdot 2^{(k+k_0)a_n}$ hold, and 
consequently 
(since $|\xi _n|^{1/a_n}\le |\xi |_a\le |\xi'|_{a'}+|\xi_n|^{1/a_n}$) 
the number $k_0$ may be taken so large that for all $k$,
\begin{equation}
  \supp \cf g_k  \subset 
   \{\,\xi\in \rn\mid
  2^{k+k_0}\le |\xi |_a\le \tfrac{11}{10} 2^{k+k_0}\,\}
  \subset \{\,\xi\mid \ffi_{k+k_0}(\xi)=1\,\}.  
\end{equation}
By the definition, the $\bspq$-norm of $g_k$ is therefore equal to
$ 2^{(k+k_0)(\frac{|a|}{p}-|a'|)}  \|\, g_k\, | L_p(\rn )\|$, so since
$\sum_{j=0}^k \ffij=\psi(2^{-ka}\cdot)$, we conclude
\begin{equation}
   \|\,g_k\, | L_p\| 
    \le c 2^{-k\frac{a_n}{p}}\|\, \cf _1^{-1}\eta \, | L_p(\re) \|  \cdot 
    2^{k(|a'|-\frac{|a'|}{p})} \|\, \cf _{n-1}^{-1}[\psi(\cdot,0)]\, | L_p(\rl )\|
\end{equation}
and the claim in \eqref{717-1} follows immediately.

To prove \eqref{717-2}, note that the definition of the Besov norm implies
\begin{equation}
\begin{split}
\label{717-11}
 \|\,g_k(\cdot ,0)\, | B^{|a'|(\frac{1}{r}-1),a'}_{r,u}(\rl)\|
  &\ge (\sum _{j=0}^{k-1} 2^{j|a'|(\frac{1}{r}-1)u}
       \|\,\cf _{n-1}^{-1}(\ffij(\cdot,0)) \, | L_r(\rl)\|^u )^{1/u}
\\
  & \ge c k^{1/u}.
\end{split}  
\end{equation}

If one assumes \eqref{5-2-1}, now in dimension $n$ again, it then follows that
$
 B^{\a(\frac{1}{p}-1),a}_{p,q}\imb A^{\a(\frac{1}{p}-1),a}_{p,q}
 \imb  B^{\a(\frac{1}{r}-1),a}_{r,\infty}
$
and this embedding between the Besov spaces implies $r\ge p$ (as one may
show by considering the functions $\rho_k$ in \cite[Lem.~4.1]{Jo1} for $-k\in\n$;
cf.\ Remark~4.6 there).
When $p\le q$ it is enough to prove $q\le r$ in the case where $r<1$, for
$q\le 1$ is a general assumption in the theorem. However, for $r<1$ there is
an inclusion $B^{\a(\frac1r-1),a}_{r,\infty}\imb
F^{0,a}_{r,\infty}$, so we get from \eqref{5-2-1} that
$A^{\a(\ep-1),a}_{p,q}\imb
F^{0,a}_{r,\infty}$ and this already implies $r\ge q$; this last
conclusion may be proved as in \cite[Prop.~3.2]{Jo3}, where $a=(1,\dots,1)$
was treated.

Step~3.
If there is an embedding as in \eqref{5-2-3}, then a Sobolev embedding gives
for some finite $t>r$ that $A^{\a(\ep-1),a}_{p,q}\hookrightarrow  
B^{\a(\frac1t-1),a}_{t,t}$, and 
this would contradict the $u$-part of \eqref{5-2-2}.
 
\subsubsection*{Proof of Proposition~\ref{ABF-prop}}
Suppose for $s:=\a(\ep-1)$ that $\aspq(\Rn)=X$, where $X$ denotes either a
Besov or a Lizorkin--Triebel space with some parameter $(\tau,\rho,\omega)$.
 
Since $\bspq\imb \aspq \imb X$ it follows from the necessity of the
parameter restrictions in the usual embeddings in the Besov and
Lizorkin--Triebel scales that 
\begin{equation}
  \tau\le s,\qquad \rho\ge p  ,\qquad
  \tau-\tfrac{\a}{\rho}\le s-\tfrac{\a}{p}.
\end{equation}
Moreover, since $X\imb B^{\a(\frac1r-1),a}_{r,\infty}$ for $r=\max(p,q)$, it
follows in the same way, since $B^{\a(\frac1r-1),a}_{r,\infty}$ has
differential dimension $-\a$, that 
\begin{equation}
  \tau-\tfrac{\a}{\rho}\ge -\a = s-\ap.
\end{equation}
Hence $\tau-\tfrac{\a}{\rho} = s-\ap$, and therefore $X$ is a Besov space
according to \eqref{5-2-3}.
 
Using these conclusions, we obtain from the necessity of \eqref{5-2-2} that
$\rho\ge r$, and because
\begin{equation}
  B^{\a(\frac1\rho-1),a}_{\rho,\omega}=
  X\imb   B^{\a(\frac1r-1),a}_{r,\infty},
\end{equation}
we find that $\rho=r$. Then we see from this and 
\eqref{5-2-1}--\eqref{5-2-2} that $\omega=\infty$. 

Finally we conclude that $\aspq=X=B^{\a(\frac1r-1),a}_{r,\infty}$ must hold
under the assumption made at the beginning of the proof; but this conclusion
is absurd since the Dirac measure $\delta_0$ belongs to
$B^{\a(\frac1r-1),a}_{r,\infty}$, which then contradicts the inclusion 
$\aspq\subset L_1$ that one obtains from Corollary~\ref{aspq-cor}.

\begin{rem}   \label{ru-rem}
The necessity of the first inequality in \eqref{5-2-6} may be obtained by
means of the special Schwartz functions $\rho_k$ constructed in
\cite[Lem.~4.1]{Jo1}; these may be inserted for $-k\in\N$ into the
inequality expressing the boundedness of the embeddings
$B^{\a(\frac1r-1),a}_{r,u}\imb A^{\a(\ep-1),a}_{p,q}\imb L_p$.
However, the second part of \eqref{5-2-6} is not so easy to handle, for when
one analogously inserts the $\rho_{N,l}^{(s-\ap)}$, the $\aspq$-norms of these
functions are troublesome to calculate.

 To prove the necessity of \eqref{5-2-8},
given \eqref{5-2-7}, it may be used for any $t<r$ that 
\begin{equation}
  B^{\a(\frac1t-1),a}_{t,t}\imb F^{\a(\frac1r-1),a}_{r,u}
  \imb A^{\a(\ep-1),a}_{p,q};
\end{equation}
from (the established necessity of) \eqref{5-2-6} it follows that $t\le
p$\,---\,hence by taking the supremum over such $t$ that $r\le p$, which is
the $r$-part of \eqref{5-2-8}. 

Clearly this would give $r\le\min(p,q)$ if
\eqref{5-2-6} could be shown to be necessary in its entirety.
So, if $r<p$ we would have deduced that $r\le q$.
Moreover, for $r=p$ the conclusion $r\le\min(p,q)$ would reduce
to the inequality $p\le q$, and so it would remain to be
proved that $u\le q$. Here one could try to calculate the norms of the
$\theta_N^{(s)}$ of \cite[Lem.~4.1]{Jo1}, but then the same difficulties
would occur as above for \eqref{5-2-6}.
\end{rem}


\appendix 
\section{Notation} \label{notation-app}
Let $\cs (\rn )$ be the Schwartz space of all complex-valued
rapidly decreasing $C^\infty$-functions on $\rn$, equipped with the usual
topology; $\cs '(\rn )$ denotes the topological dual, the space
of all tempered distributions on $\rn$.
If $\ffi\in \cs (\rn )$ then $\widehat {\ffi}= \cf \ffi$ and
$\check{\ffi}=\cf ^{-1}\ffi $ are respectively the Fourier and the inverse
Fourier transform of $\ffi$, extended in the usual way 
from $\cs (\rn )$ to $\cs '(\rn )$. 

The space of uniformly continuous, bounded functions on $\rn$, valued in
a Banach space $X$, is denoted by
$C_{\op{b}}(\rn, X)$; for $X=\C$ the Banach space is suppressed.

For a normed or quasi-normed space $X$ we denote by $\|x\,|\,X\|$ the norm of 
the vector $x$. Recall that $X$ is quasi-normed when the triangle inequality
is weakened to $\|x+y\,|\,X\|\le c(\|x\,|\,X\|+\|y\,|\,X\|)$ for some
$c\ge 1$ independent of $x$ and $y$.

\medskip

All unimportant positive constants are denoted by $c$, occasionally with
additional subscripts within the same formulas. 
The equivalence ``$\mbox{term}_1 \sim\mbox{term}_2$'' means that there exist
two constants $c_1, c_2 >0$ independent of the variables in the two terms
 such that $c_1\, \mbox{term}_1\le \mbox{term}_2\le c_2\, \mbox{term}_1$.

\section{Anisotropic function spaces}   \label{spaces-app} 
The conventions we adopt here are, by and large, those of \cite{Ya}. 
For each coordinate $x_i$ in $\rn$  a weight
$a_i$ is given such that   $\min(a_1,\dots,a_n)=1$.
The vector $a=(a_1,...,a_n)$  is called an $n$-dimensional anisotropy,
and $|a|:=a_1+...+a_n$.
If $a=(1,...,1)$ we have  the {\it isotropic}  case.

For a given $a=(a_1,...,a_n)$, the action of $t\in [0,\infty )$ on $x\in \rn$
is defined by the formula
\begin{equation}
  \label{112-2}
  t^ax=\left (t^{a_1}x_1,...,t^{a_n}x_n \right ).
\end{equation}
For $t>0$ and $s\in \re$ we set $t^{sa}x=(t^{s})^{a}x$. In particular,
$t^{-a}x=(t^{-1})^{a}x$ and $2^{-ja}x=(2^{-j})^{a}x$. 

For $x=(x_1,...,x_n)\in \rn$, $x\neq 0$, let
 $|x|_a$ be the unique positive number $t$ such that
\begin{equation}
\label{112-6}
\frac{x_1^2}{t^{2a_1}}+...+\frac{x_n^2}{t^{2a_n}}=1
\end{equation}
and let $|0|_a=0$ for $x=0$.

By \cite[1.4/3,8]{Ya},  the map
$|\cdot|_a$ is an anisotropic distance function, which is $C^\infty$ and
coincides with $|\cdot|$ in the isotropic case. 
(Anisotropic distance functions are continuous
 maps $u \colon\rn\to \re$  fulfilling
 $u(x)>0$ if $x\neq 0$ and 
 $u(t^{a}x)=t u(x)$ for all $t>0$ and all $x\in\rn$;
 any two such functions $u$ and $u'$ are equivalent
 in the sense that $u(x)\sim u'(x)$, see \cite{st-wai} and \cite[1.2.3]{di}.) 
Moreover, setting $a_{max}=\max\{\,a_i\mid 1\le i\le n\,\}$ one has, cf.\
\cite[1.4]{Ya}, for any $x\in\rn$ that
\begin{equation}
  \min \{\,|x|, |x |^{1/a_{max}}\,\}\le
  |x |_a \le
  \max \{\,|x |, |x |^{1/a_{max}}\,\}.
\end{equation}

\medskip

If $(\varphi_j)_{j\in\n}$ is the anisotropic partition of unity from Section~\ref{gen-sect},
then for any $f\in \cs '(\rn )$,
\begin{equation}
\label{700-0}
f=\jsum \cf ^{-1}( \ffij \cf f) \qquad \mbox{with convergence in }
\quad \cs ' (\rn ).
\end{equation}
Let $0<p\le \infty$, $0<q\le \infty$, $s\in \re$. The anisotropic Besov space 
$\bspq (\rn )$ and, provided $p<\infty$, the anisotropic Lizorkin--Triebel
space $\fspq (\rn )$ are defined to consist of all tempered distributions
$f\in \cs '(\rn )$ for which the following quasi-norms are finite:
\begin{gather}
\label{211-1}
 \|f\, |\, \bspq (\rn )\|=
                     (\jsum 2^{jsq}
                \|\cfi (\ffij  \cf {f} )\, |  L_{p}(\rn )\|^{q}
                  \ )^{1/q},
\\
\label{211-2}
\|f\,|\,\fspq (\rn )\|=
 \bigl\| 
 (\jsum 2^{jsq}|\cfi (\ffij \cf {f}) (\cdot )|^{q} )^{1/q}
  \, \bigm|  L_p (\rn ) \bigr\|,
\end{gather} 
respectively (with the usual modification if $q=\infty$).

Both $\bspqa (\rn )$ and $\fspqa (\rn )$ are quasi-Banach spaces (Banach spaces if 
$p\ge 1$ and $q\ge 1$) which  are independent of the choice of  $(\ffij )_{j\in \no}$.        

The embeddings $\cs (\rn )\imb\bspqa (\rn )\imb \cs '(\rn )$ and 
$\cs (\rn )\imb\fspqa (\rn )\imb \cs '(\rn )$ hold true for all admissible 
values of $p,q,s$. Furthermore, if both
$p<\infty$ and $q<\infty$, the ranges of these inclusions are all dense;
cf.~\cite[3.5]{Ya} and \cite[1.2.10]{di}. 
The results on embeddings are reviewed (and extended) in
Appendix~\ref{embd-app} below; let us conclude with a few identifications.

If $1<p<\infty$ and $s\in\re$ then 
$F^{s,a}_{p,2}(\Rn )=H^{s,a}_p (\Rn )$ with equivalent quasi-norms; hereby 
\begin{equation}
H^{s,a}_{p}(\Rn )=
\bigl \{\,
  f\in \cs '(\rn )\bigm|
  \|\,(\sum_{k=1}^{n}(1+\xi _k^{2})^{s/(2a_k)} \widehat{f}
    )^{\vee} \, | L_p(\rn ) \| <\infty
\,\bigr\}  
\end{equation}
is the anisotropic Bessel potential space; cf.\ \cite[Rem.~11]{sto-tri}, 
\cite[2.5.2]{tri-f}, and  \cite[3.11]{Ya}.

Furthermore, if $1<p<\infty$, $s\in\re$ and if  
$s_1=\frac{s}{a_{1}}\in \N$,..., $s_n=\frac{s}{a_n}\in \N$ then 
$F^{s,a}_{p,2}(\rn )=W_p^{s,a}(\rn )$ (with equivalent quasi-norms), where
\begin{equation}
W_p^{s,a}(\Rn )=
\bigl \{\, f\in \cs '(\rn )\bigm|
  \|\,f\, |L_p (\rn )\| + \sum_{k=1}^{n}
   \|\, \frac{\partial ^{s_k}f}{\partial x_k^{s_k}}\, |L_p (\rn )
  \| <\infty
\,\bigr\}  
\end{equation}
is the classical anisotropic Sobolev space on $\rn$. 
If  $s>0$ and $\frac{s}{a_k}\notin \N$ for $k=1,...,n$,
 then  $B_{\infty,\infty}^{s,a}(\rn )= C^{s,a}(\rn )$ are the 
anisotropic H\"older spaces.

Anisotropic function spaces have been intensively studied
by { S. M. Nikol'skij}, \cite{Nik}, and by 
{ O. V. Besov}, { V. P. Il'in} and { S. M. Nikol'skij}, \cite{BIN}.
See also works of
M. Yamazaki \cite{Ya}, H.-J. Schmeisser and H. Triebel
\cite[4.2]{ST},   
A. Seeger \cite{see}, P. Dintelmann \cite{di,din}
etc.

\section{Properties of the anisotropic spaces}
  \label{atom-app}
The point of this section is to sketch how Propositions~\ref{711-t21} and
\ref{indep-prop} may be proved in the present anisotropic set-up.
Along the way, we also prove some interpolation formulas that are of interest
in their own right.
The main tool will be anisotropic versions of the atomic decompositions of
\cite{FJ1,FJ2}, see  W.~Farkas \cite{Far1}.

\subsection{Atomic decompositions}
As a preparation, we recall some basic notions
of atomic decompositions in an anisotropic setting. 

Consider the lattice $\zn$ as a subset of $\rn$.
If $\nu\in\no $ and $m=(m_1,...,m_n)\in \zn$,
we denote by $\qnium$ the rectangle in $\rn$ centred at $2^{-\nu a}m=
(2^{-\nu a_1}m_1,...,2^{-\nu a_n}m_n)$, 
which has sides parallel to the axes and side lengths
$2^{-\nu a_1}$, ..., $2^{-\nu a_n}$, respectively ($Q_{0m}^{a}$ is a cube
with side length 1).
If $\qnium$ is such a rectangle in $\rn$
and $c>0$, then $c\qnium$ denotes the concentric rectangle with side lengths
$c 2^{-\nu a_1}$, ..., $c 2^{-\nu a_n}$.
If $\beta =(\beta _1,...,\beta _n)\in \non$ is a multi-index 
 and if $x=(x_1,...,x_n)\in \rn$, then 
$x^{\beta}:=x_1^{\beta _1}\cdots x_n^{\beta _n}$, and 
we write $a\beta =a_1\beta _1+...+a_n\beta _n$.
If $E\subset\rn$ is Lebesgue measurable, then
$|E|$ denotes its Lebesgue measure.
Now, anisotropic atoms are defined as follows.

\begin{defi}
\label{def-21}
 Let $s\in \re$, $0<p\le\infty$ and $K$, $L\in \re$. A function
 $\rho\colon\rn\to \C$ for which $D^{\beta}\rho $ exists
 when $a\beta \le K$ (or $\rho$ is continuous if $K\le 0$) 
is called an anisotropic $(s, p)_{K,L}$-atom, if
\begin{gather}
   \supp\rho \subset c\qnium
   \quad\text{for some $\nu\in \N$, $ m\in\zn$ and $c>1$} ,
    \label{312-2a}
\\  
 | D^{\beta}\rho (x)|^{\a}\le |\qnium |^{s-a\beta-\frac{\a}{p}}
 \quad\mbox{if} \quad a\beta\le K,
  \label{312-2b}
\\
  \label{312-2c}
 \int _{\rn}x^{\beta}\rho (x) dx=0\quad \mbox{if} \quad a\beta \le L.
\end{gather}
If conditions {\em (\ref{312-2a})} and {\em (\ref{312-2b})}
are satisfied for $\nu =0$, then  $\rho$
is called an anisotropic $1_K$-atom.
\end{defi} 

If the atom $\rho$ is located at $Q_{\nu m}^a$ 
(i.e.\ $\supp \rho\subset c \qnium$ with $\nu\in \no\, $,
$m\in \zn$, $c>1$) we denote it by $\rho_{\nu m}^{a}$.
The value of the number $c>1$ in \eqref{312-2a} is unimportant; 
it allows, at any level $\nu$, a controlled overlap of the supports of
different $\rho _{\nu m}^{a}$.

The main advantage of the atomic approach is that one 
can often reduce a problem  given in $\bspq$ or $\fspq$
to some corresponding sequence spaces; these are here denoted by
$b_{p,q}$ and $f^{a}_{p,q}$.
If $\qnium$ is a rectangle as above, let
$\chi _{\nu m}$ be the characteristic function of $\qnium$; then
\begin{equation}
  2^{\nu |a|/p}\chi _{\nu m}
\end{equation}
is the $L_p(\rn )$-normalised characteristic function of $\qnium $ whenever
$0<p\le\infty$. 

If $0<p,q \le\infty$, then $b_{p,q}$ is the collection of all 
$\lambda =\{\,\lambda _{\nu m}\in \C \mid \nu\in\no\, , m\in\zn\,\}$
such that
\begin{equation}
\label{312-6}
  \|\,\lambda \,|b_{p,q}\| =
 \bigl ( \niusum  ( \msum |\lambda _{\nu m}|^{p} )^{q/p} \bigr )^{1/q}   
 \end{equation}
is finite (usual modification if $p=\infty$ and/or $q=\infty$).
Furthermore, $f_{p,q}^{a}$ is the collection of all such sequences 
$\lambda$ for which
\begin{equation}
\label{312-7}
  \|\,\lambda \,| f_{p,q}^{a}\| =
 \bigl \| ( \niusum \msum |\lambda _{\nu m}
 2^{\nu |a|/p}\chi _{\nu m}(\cdot )|^{q})^{1/q} \, \bigm| L_p (\rn ) \bigr \|
\end{equation}
(usual modification if $p=\infty$ and/or $q=\infty$) is finite.

For $0<p\le\infty$ and $0<q\le\infty$ we will use the abbreviations
\begin{equation}
\label{31*-0}
\sigma _p=|a|\bigl (\frac{1}{p}-1 \bigr)_{+}\quad\mbox{and}\quad
\sigma _{p,q}=|a|\bigl (\frac{1}{\mbox{min}(p,q)}-1\bigr )_{+}.
\end{equation}

\begin{prop} \label{theo-23}
{\rm (i)} Let $0<p<\infty$, $0<q\le \infty$ and $s\in \re$,
and let $K$, $L\in \re$ fulfil
\begin{gather}
\label{33-2}
 K\ge a_{max}+s\quad\mbox{if}\quad s\ge 0,
\\
\label{33-3}
L\ge \sigma _{p,q}-s.
\end{gather}
Then $g\in \cs '(\rn )$ belongs to $\fspqa (\rn )$ if, and only if, 
for some $\lambda=(\lambda_{\nu m})$ in $f^a_{p,q}$,
\begin{equation}
\label{33-7}
g=\niusum\msum\lambda _{\nu m}\rho _{\nu m}^{a}, \
\text{with convergence in $\cs '(\rn )$},
\end{equation}
where $\rho _{\nu m}^{a}$ are anisotropic $1_K$-atoms ($\nu =0$) or
anisotropic $(s,p)_{K,L}$-atoms ($\nu\in\N$).

Furthermore, $\inf \|\,\lambda\,|f_{p,q}^{a}\|$ 
with the infimum taken over all admissible representations 
\eqref{33-7} is an equivalent quasi-norm in $\fspq (\rn )$.

{\rm (ii)\/} The analogous statements are valid for the Besov spaces
$\bspq(\Rn)$ for $0<p\le\infty$ provided $\sigma_{p,q}$ and $f^a_{p,q}$
(together with its norm $\|\,\cdot\,| f^{a}_{p,q}\|$) are replaced by
$\sigma_p$ and $b_{p,q}$, respectively.
\end{prop}

\noindent
{\bf  Proof } 
The proposition is a slightly different version of the atomic decomposition 
theorem proved in  \cite{Far1}; the modifications needed are immaterial, so 
we omit details.

\subsection{Real Interpolation}

Our aim is to prove a refined Sobolev embedding due
to B.~Jawerth \cite{Ja} and J.~Franke \cite{Fr} in the isotropic context. 
As usual $(\cdot ,\cdot )_{\theta, q}$ denotes the real interpolation.

\begin{lem}
Let $s_0\neq s_1$
and $s=(1-\theta )s_0+\theta s_1$, where $0<\theta <1$.
If $0<p\le \infty$ then
\begin{gather}
\label{700-31}
(B^{s_0,a}_{p,q_0}(\rn ), B_{p,q_1}^{s_1,a}(\rn ))_{\theta ,q}=
\bspq (\rn ),
\\
\label{700-32}
(F^{s_0,a}_{p,q_0}(\rn ),F_{p,q_1}^{s_1,a}(\rn ) )_{\theta ,q}
=\bspq (\rn ),
\end{gather}
provided $0<p<\infty$ in the last formula.
\end{lem}

Formulas  (\ref{700-31}) and (\ref{700-32})
  can be proved using the same arguments as in
\cite[2.4.2]{Tr6}, and we refrain from doing this here.
In the case of constant $s$, there is another result:

\begin{prop}
\label{711-t23}
Let $0<\theta <1$ and $0<p_0<p<p_1<\infty$. When
$
 \frac{1}{p}=\frac{1-\theta}{p_0}+\frac{\theta}{p_1}, 
$
then 
\begin{equation}
\label{700-33}
(F_{p_0,q}^{s,a}(\rn ), F_{p_1,q}^{s,a}(\rn ))_{\theta ,p}=
\fspq (\rn ).
\end{equation}
\end{prop}

\noindent
{\bf Proof }
Replacing (if necessary) $|\, \cdot \, |_a$ by an
equivalent anisotropic distance function we may assume
\begin{equation}
  \{\,x\in \rn \mid |x|_a\le 2\,\}\subset [-\pi ,\pi ]^{n}.  
\end{equation}
Let $f_{p,q}^a$ as in (\ref{312-7}).
By immaterial modifications of the proof of \cite[(6.10)]{FJ2} 
we have
\begin{equation}
\label{700-050}
(f_{p_0,q}^a, f_{p_1,q}^a )_{\theta ,p}=f_{p,q}^a \quad \mbox{for}\quad
\frac{1}{p}=\frac{1-\theta}{p_0}+\frac{\theta}{p_1}.
\end{equation}
Then \eqref{700-33} is  a  consequence
of \eqref{700-050} and the following anisotropic $\varphi$-transform.

\begin{prop}
\label{ap2-theo-44}
Let $(\ffij )_{j\in \no}$ be a partition of unity as in Section~\ref{gen-sect},
 and let $\rho \in \cs(\rn)$ fulfil $\rho (x)=1$ if $|x|_a\le 2$ and
$\supp\rho\subset [-\pi ,\pi ]^{n}$. 
The operators $U_{\ffi} :\fspq (\rn )\to f_{p,q}^{a}$
and $T_{\rho} :  f_{p,q}^{a}\to\fspq (\rn )$ defined by
\begin{equation}
  \label{ap2-342-2}
  U_{\ffi}(g) = \bigl \{\,
  (2\pi )^{-n/2}\, 2^{\nu\left (s-\frac{n}{p}\right )} (\ffi_{\nu}\widehat{g})^{\vee}
  (2^{-\nu a}m)\bigm|  {\nu\in \no\, , m\in\zn}
  \,\bigr\}
\end{equation}
for $g\in\fspq (\rn )$ and by 
\begin{equation}
  \label{ap2-342-3}
  T_{\rho} (\lambda ) = \niusum \,\msum \lambda _{\nu m}\, 2^{-\nu \left (s-\frac{n}{p}\right )}
  \check{\rho}(2^{\nu a}\cdot -m)
\end{equation}
for $\lambda =\{\, \lambda _{\nu m}\mid {\nu\in \no \, , m\in\zn}\,\} $
belonging to $f^a_{p,q}$, respectively,
are bounded.

Furthermore, 
$(T_{\rho}\, \circ\, U_{\ffi})(g)=g$ for any $g\in \fspq (\rn )$ and 
$\|U_{\ffi}(\cdot )\,|\,f_{p,q}^{a}\|$
is an equivalent quasi-norm on $\fspq (\rn )$.

The corresponding results hold for the $\bspq (\rn )$ 
spaces with $b_{p,q}$ in place of $f_{p,q}^{a}$.
\end{prop}

\noindent
{\bf Proof } For $p,q<\infty$ a proof may be found in \cite{din}, where
density of $\cs (\rn )$ in $\bspq (\rn )$ and $\fspq (\rn )$ was used.
In the remaining cases one may proceed, for example, 
as in \cite[14.15]{tri97}.

The proposition, due to P.~Dintelmann \cite[Theorem 1]{din},
represents
the anisotropic version of a theorem originally proved by
M.~Frazier and B.~Jawerth, see \cite{FJ2}. 

\subsection{Embeddings} 
  \label{embd-app}
In addition to elementary embeddings (monotonicity in $s$ and $q$) we have 
\begin{equation}
\label{700-8}
B_{p, \min (p,q)}^{s,a}(\rn )\imb \fspq (\rn )
\imb B_{p, \max (p,q)}^{s,a}(\rn );
\end{equation}
hence $B^{s,a}_{p,p}=F^{s,a}_{p,p}$ whenever $0<p<\infty$.
There are, moreover, the Sobolev embeddings 
\begin{equation}
\label{700-21}
\bspq (\rn )\imb B_{r,q}^{t,a}(\rn )
\qquad \mbox{and} \qquad
       \fspq (\rn )\imb F_{r,u}^{t,a}(\rn )
\end{equation}
provided that
\begin{equation}
p < r \qquad \mbox{and} \qquad
        s-\frac{|a|}{p} = t-\frac{|a|}{r};  
\end{equation}
$q$ and $u$ are independent of each other. These assertions may
conveniently be found in \cite{Ya}.
There is also an anisotropic version of the
Jawerth--Franke embedding (which has features in common with both of the
above types):

\begin{prop}
\label{711-t21}
Let $0<p_0<p<p_1\le \infty$, $s_1<s<s_0$ and $0<q\le \infty$. Then 
\begin{equation}
\label{700-23}
B_{p_0,p}^{s_0,a}(\rn )\imb 
F_{p,q}^{s,a}(\rn )
\imb B_{p_1,p}^{s_1,a}(\rn )
\end{equation}
provided $s_0-\frac{|a|}{p_0} =s-\frac{|a|}{p}=s_1-\frac{|a|}{p_1}$.
\end{prop}

\noindent
{\bf Proof }
Let $0<p_0<p'<p<p''<\infty$ and let
\begin{equation}
s'=s+\frac{|a|}{p_0}-\frac{|a|}{p'}   \quad \mbox{and}\quad 
  s''=s+\frac{|a|}{p_0}-\frac{|a|}{p''}.   
\end{equation}
As a consequence of (\ref{700-21}) we obtain
\begin{equation}
   F_{p_0,q}^{s',a}(\rn )\imb F_{p',q}^{s,a}(\rn )\quad
   \mbox{and}\quad
   F_{p_0,q}^{s'',a}(\rn )\imb F_{p'',q}^{s,a}(\rn ).  
\end{equation}
There exists $\theta \in\,]0,1[\,$ such that 
$\frac{1}{p}=\frac{1-\theta}{p'}+\frac{\theta}{p''}$.
Then $s_0=(1-\theta )s'+\theta s''$.

Using (\ref{700-32}), elementary properties of the interpolation,
and (\ref{700-33}) we have
\begin{equation}
  B_{p_0,p}^{s_0,a}(\rn )=
  (F_{p_0,q}^{s',a}(\rn ), F_{p_0,q}^{s'',a}(\rn ))_{\theta ,p}\imb
  (F_{p',q}^{s,a}(\rn ), F_{p'',q}^{s,a}(\rn ))_{\theta , p}= 
  \fspq (\rn )  
\end{equation}
and this gives the first embedding in \eqref{700-23}.

To prove the second, let now $p'<p<p''<p_1\le\infty$ and
$s^{(j)}=s-\frac{|a|}{p^{(j)}}+\frac{|a|}{p_1}$ for $j=1,~2$.
Then by (\ref{700-21}) and (\ref{700-8}) we have
$F^{s,a}_{p^{(j)},q}(\rn )\imb B_{p_1,\infty}^{s^{(j)},a}(\rn )$, and
$s_1=(1-\theta )s'+\theta s''$ when $\frac{1-\theta}{p'}+\frac{\theta}{p''}=\frac{1}{p}$,
so applying \eqref{700-33} and \eqref{700-31} we conclude that
\begin{equation}
F_{p,q}^{s,a}(\rn )=
  (F^{s,a}_{p',q}(\rn ), F^{s,a}_{p'',q}(\rn ))_{\theta ,p} \imb
 (B_{p_1,\infty}^{s',a}(\rn ), B_{p_1,\infty}^{s'',a}(\rn ))_{\theta ,p}
 = B_{p_1,p}^{s_1,a}(\rn ).  
\end{equation}

\begin{rem}
The second and first part of \eqref{700-23} was proved by B.Jawerth \cite{Ja}
and by J. Franke \cite{Fr}, respectively, in the isotropic case. The present
extension to $a\neq (1,\dots,1)$  would 
have been beneficial for e.g.\ the fine continuity properties of the pointwise
product as investigated in \cite{Jo1}, where it was necessary to distinguish
between the isotropic and anisotropic cases in numerous places.
\end{rem}

\subsection{The $q$-independence of the range space $\gamma_0 (\fspq(\rn))$}
  \label{qindep-app}

\begin{prop}   \label{indep-prop}
When $\g$ is defined on both $\fspq(\Rn)$ and $F^{s,a}_{p,t}(\Rn)$, then
\begin{equation}
\g (\fspq(\Rn) ) =  \g (F^{s,a}_{p,t}(\Rn )).  
\end{equation}
\end{prop}

\noindent
{\bf Proof } 
We proceed as in \cite[11.1]{FJ2}.
If $q<t$, the elementary embedding $\fspq(\Rn) \imb F^{s,a}_{p,t}(\Rn )$  implies
$ \g ( \fspq(\Rn) )\subset  \g ( F^{s,a}_{p,t}(\Rn ))$.
To prove the converse inclusion, let $K$ and $L$ fulfil \eqref{33-2}--\eqref{33-3}
and let us write $g\in F^{s,a}_{p,t}(\rn )$ as
\begin{equation}
g=\niusum\,\msum \lambda _{\nu m} \rho _{\nu m}^{a,t}\quad ,
\quad \mbox{convergence  in}\quad \cs '(\rn ),  
\end{equation}
with
$\|\lambda \,|\,f_{p,t}^{a}\|\le c\;\|g\,|\,F^{s,a}_{p,t}(\rn )\|$;
cf.~\eqref{33-7}.
We claim there exists a function $\widetilde{g}\in \fspq(\Rn)$ with 
$ \g (g) = \g (\widetilde{g})$.
Only such rectangles $\qnium$ for which the relevant $c\qnium$ 
intersects $\Gamma=\{\, x\mid x_n=0\,\}$, are important.
Let
\begin{equation}
\label{z32-4}
A=\bigl\{\,(\nu ,m)\in \no\times\zn\bigm| c\qnium\cap \Gamma
      \neq \emptyset\,\bigr\}.
\end{equation}
We define 
$\widetilde{\lambda}_{\nu m}=\lambda _{\nu m}$ if $(\nu ,m)\in A$ and 
otherwise $\widetilde{\lambda}_{\nu m}=0$, and put 
$\widetilde{\lambda}=\{\,\widetilde{\lambda}_{\nu m}\mid
{\nu\in \no\, ,m\in \zn}\,\}$.
Let now $\psi\in S(\re)$ be such that  $\supp\psi \subset 
\left [-\tfrac{1}{2},\tfrac{1}{2}\right ]$ and $\psi (0)=1$ and, moreover, 
\begin{equation}
  \int _{\re}z^{\beta _n}\psi (z) dz =0\quad\mbox{for all}\quad
  \beta _n\in \no\quad\mbox{such that}\quad a_n\beta _n\le L.  
\end{equation}
Using this we define
\begin{equation}
   \widetilde{\rho}_{\nu m}^{a,q}(x',x_n)={\rho}_{\nu m}^{a,t}(x',0)\,\psi
 (2^{\nu a_n}x_n)
\end{equation}
and remark that $\widetilde{\rho}_{\nu m}^{a,q}$ is supported in a rectangle 
$c\widetilde{Q}_{\nu m}^{a}$ where $\widetilde{Q}_{\nu m}^{a}$ 
has sides parallel to the axes,
is centred at $(2^{-\nu a_1}m_1,...,2^{-\nu a_{n-1}}m_{n-1}, 0)$ and
its side lengths are respectively $2^{-\nu a_1}$,...,$2^{-\nu a_{n-1}}$,
 $2^{-\nu a_n}$.
Furthermore, if $\beta =(\beta ',\beta _n)\in \non$
is such that $a\beta \le K$ and if $\nu\in \N$ then
\begin{equation}
|D^{\beta }\widetilde{\rho}_{\nu m}^{a,q}(x)|\le 
c\;|D^{\beta '}\rho _{\nu m}^{a,t}(x',0)|\cdot
        2^{\nu a_n\beta _n}
        \le c'\;2^{-\nu\left (s-\frac{|a|}{p}\right )}\,
         2^{\nu a\beta}.  
\end{equation}
It follows that each $\widetilde{\rho }_{\nu m}^{a,q}$,
up to an unimportant constant,
is an anisotropic $1_K$-atom for $\nu=0$ or an anisotropic $(s,p)_{K,L}$-atom
(due to its product structure and to the assumptions on the function
$\psi$ there are no problems in checking the moment conditions).

Defining
\begin{equation}
  \label{z32-24}
  \widetilde{g}=\niusum\, \msum\widetilde{\lambda}_{\nu m}\widetilde{\rho}_{\nu m}^{a,q}
  =\sum \limits _{(\nu ,m)\in A}{\lambda}_{\nu m}\widetilde{\rho}_{\nu m}^{a,q}
\end{equation}
we have $\g (g) = \g (\widetilde{g})$. For $(\nu , m)\in A$, let 
\begin{equation}
 \widetilde{E}_{\nu m}^{a}=
  \bigl\{\,(x_1,...,x_n)\in \widetilde{Q}_{\nu m}^{a}\bigm| 
   2^{-(\nu +1)a_n-1}<x_n\le  2^{-\nu a_n-1} \bigr\}.  
\end{equation}
Obviously, 
\begin{equation}
  \frac{|\widetilde{E}_{\nu m}^{a}|}{|\widetilde{Q}_{\nu m}^{a}|}= 
  \frac{1-2^{-a_n}}{2}>0.  
\end{equation}
Then (with the usual modification for $q=\infty$)
\begin{equation}
\label{z32-25}
\|\widetilde{\lambda}\,|\,f_{p,q}^{a}\|
\sim \bigl\|
        \bigl (\sum\limits _{(\nu , m)\in A}|\lambda _{\nu m}\,
        \widetilde{\chi}_{\nu m}^{(p)}
           (\cdot )|^{q}\bigr )^{1/q}\,\bigm| L_p (\rn )
     \bigr \|
\end{equation}
where $\widetilde{\chi}_{\nu m}^{(p)}$
denotes the $L_p(\rn )$-normalised characteristic function of 
the rectangle $\widetilde{E}_{\nu m}^{a}$.
Using an inequality of Fefferman--Stein type, see \cite{FS} for the anisotropic
Hardy--Littlewood maximal function, the proof is a simple anisotropic 
counterpart of  \cite[2.7]{FJ2}. 

For $(\nu ,m)\in A$ the sets $\widetilde{E}_{\nu m}^{a}$ 
are pairwise disjoint and so at most
one term in the sum on the right-hand side is nonzero. 
Hence $q$ and $1/q$ cancel in \eqref{z32-25} and may therefore be replaced by $t$ and $1/t$. 
So, with the usual modification if $t=\infty$,
\begin{equation}
\|\widetilde{\lambda}\,|\,f_{p,q}^{a}\|\sim 
\bigl\| \bigl ( \sum _{(\nu , m)\in A}|\lambda _{\nu m}\,
 \widetilde{\chi}_{\nu m}^{(p)}
           (\cdot )|^{t}\bigr )^{1/t}\,\bigm|L_p (\rn )\bigr\|\le c\;
\|{\lambda}\,|\, f_{p,t}^{a}\|.
\end{equation}
The last relation together with 
(\ref{z32-24}) prove the fact that $\widetilde{g}\in \fspq (\rn )$ and 
\begin{equation}
\label{qt-est}
\|\widetilde{g}\,|\, \fspq (\rn )\|\le c\;
\|\widetilde{\lambda}\,|\,f_{p,q}^{a}\|\le
c'\;\|{\lambda} \,|\, f_{p,t}^{a}\|\le c''\;\|g\,|\,F^{s,a}_{p,t}(\rn )\|.
\end{equation}
This verifies our claim.

\begin{rem}
The $q$-independence of the traces of Lizorkin--Triebel spaces has also been
treated by Yu.~V.~Netrusov~\cite{Ne2}.
\end{rem}

 \section{Series of entire functions with compact spectra}
   \label{series-app}
 First a few well-known results on convergence of certain series are recalled:

 \begin{lem}   \label{cnv1-lem}
 Let $0 < A < \infty$ and let $\{f_k\}_{k=0}^\infty $ be a sequence of
 functions on $\rn$ such that
 \begin{equation}\label{5-1-1}
 \left. \begin{array}{ll}
 &\supp \cf f_k \subset \{\, \xi \mid  \tfrac1A 2^k \le |\xi |_{a} \le A
 \, 2^k \,\} , \qquad k = 1 , \ldots ,\\
 &\supp \cf f_0 \subset \{\, \xi\mid  |\xi |_{a} \le  A
 \,\} .
 \end{array}\right\}
 \end{equation}
 Then one has, for all $s \in \re$ and $0 < p \le\infty$, that
 \begin{equation}\label{5-1-3}
 \| \, \sum_{k =0}^\infty \, f_k \, | \bspq (\Rn) \| \le c \, 
 \bigl( \sum_{k =0}^\infty \, 2^{skq} \| \, f_k  \, | L_p(\rn) \|^q
 \bigr)^{1/q} .
 \end{equation}
 More precisely, if the right-hand side is finite, then $\sum_{k=0}^\infty
 f_k$ converges in $\cs'$ to a distribution satisfying this inequality
 (where $c$ depends on $a,A,p,s,n$ but not on $\{f_k\}$).

 For $p<\infty$ an analogous result holds for $\fspq$, provided that (on the
 right hand side of \eqref{5-1-3}) the $\ell_q$-norm is
 calculated pointwise at each $x\in\rn$ before the $L_p$-norm is taken.
 \end{lem}

 \begin{lem}    \label{cnv2-lem}
 Let $0 < A < \infty$ and  let $\{f_k\}_{k =0}^\infty$ be a sequence of
functions on $\rn$ such that
\begin{equation}\label{5-1-4}
\supp \cf f_k \subset \{\, \xi \mid   |\xi |_{a} \le A
\, 2^k \,\}\, , \qquad k =0, 1, \ldots ,
\end{equation}
and suppose that
\begin{equation}\label{5-1-5}
s  > \a (\tfrac{1}{p} - 1)_+  .
\end{equation}
Then the statements in \eqref{5-1-3} ff.\ hold true. Moreover, for
$p<\infty$ one has the analogous result for the space $\fspq$ provided
\eqref{5-1-5} is replaced by the condition 
\begin{equation}\label{5-1-5'}
s  > \a\max (\tfrac{1}{p} - 1, \tfrac{1}{q}-1,0 ) .
\end{equation}
\end{lem}
These lemmas are proved in \cite{Ya}, but see also 
\cite{Jo1},\cite{Mar2} or \cite[2.3.2]{RS}.
For the borderline case with equality in \eqref{5-1-5} we refer to
\cite[Th.~3.1]{Jo3}; it is straightforward to get anisotropic variants of
this result, so without proof we state what is needed above for the case $p<1$: 

\begin{lem}  \label{cnv3-lem}
Let $\{f_k\}_{k =0}^\infty$ be a sequence of
functions such that \eqref{5-1-4} is satisfied for some $A<\infty$.
Let $0 < p < 1$, $0<q\le1$ and suppose that
$ s  = \a (\frac{1}{p} - 1) $.
Then we have 
\begin{equation}\label{5-1-12}
\| \, \sum_{k =0}^\infty f_k \, | L_1(\Rn) \| \le c \, 
 \big( \sum_{k =0}^\infty  \, 2^{k\a (\frac 1 p - 1)q} \|\, f_k
\, |   L_p(\Rn) \|^q\big)^{1/q}  .
\end{equation}
More precisely, if the right-hand side is finite, then $\sum_{k=0}^\infty
f_k$ converges in $L_1(\Rn)$ to a distribution $f$ satisfying this inequality
(where $c$ depends on $A,p,s,n$ but not on $f_k, k=0, 1 , \ldots $).
Moreover, this limit $f$ also satisfies 
\begin{equation}\label{5-1-13}
  \| \, \sum_{k =0}^\infty f_k \, | {B}^{\a (\frac 1 r - 1),a}_{r,\infty}(\rn) \|
  \le c \,  
  \big( \sum_{k =0}^\infty  \, 2^{k\a (\frac 1 p - 1)q} \| \, f_k 
  \, |  L_p(\Rn) \|^q\big)^{1/q}  
\end{equation}
for $r=\max(p,q)$.
\end{lem}

It is of course a stronger fact that the sum $f$ belongs to
the space $\aspq$, but the lemma is needed for the proof of \eqref{bab-eq} above.




\end{document}